\documentclass[a4paper,12pt]{article}

\usepackage{float}
\usepackage{amsmath} 
\usepackage{amssymb}
\usepackage{amsfonts}
\usepackage{epsfig}
\usepackage{color,fancybox,graphicx}

\numberwithin{equation}{section}

\newtheorem{rem}{Remark}[section]

\newtheorem{lem}{Lemma}[section]
\newtheorem{cor}{Corollary}[section]

\newtheorem{theo}{Theorem}[section]

\newcommand{\bX}{\mathbf{X}}
\newcommand{\bx}{\mathbf{x}}
\newcommand{\bW}{W}

\begin{document}

\begin{center}

{\sc \Large Sparse Single-Index Model\\

\vspace{0.5cm}}

Pierre Alquier\footnote{Corresponding
author.} 
\vspace{0.2cm}\\
LPMA\footnote{Research partially supported by the French ``Agence
Nationale pour la Recherche''
under grant ANR-09-BLAN-0128 ``PARCIMONIE''.}\\
Universit\'e Paris Diderot -- Paris VII\\
Bo\^{\i}te 188, 175 rue du Chevaleret\\
75013 Paris, France\\

\vspace{0.2cm}
CREST-LS\\
3 avenue Pierre Larousse\\
92240 Malakoff, France\\
\smallskip
\texttt{alquier@math.jussieu.fr}
\vspace{0.5cm}

G\'erard Biau
\vspace{0.2cm}\\
LSTA \& LPMA\footnote{Research partially supported by the French ``Agence
Nationale pour la Recherche''
under grant ANR-09-BLAN-0051-02 ``CLARA''.}\\
Universit\'e Pierre et Marie Curie -- Paris VI\\
Bo\^{\i}te 158,  Tour 15-25, 2\`eme \'etage\\
4 place Jussieu, 75252 Paris Cedex 05, France\\

\vspace{0.2cm}
DMA\footnote{Research carried out within the INRIA project ``CLASSIC'' hosted by
Ecole Normale Sup{\'e}rieure and CNRS.}\\
Ecole Normale Sup\'erieure\\
45 rue d'Ulm\\
75230 Paris Cedex 05, France\\
\smallskip
\texttt{gerard.biau@upmc.fr}
\end{center}

\vspace{0.05cm}
\begin{abstract} 
\noindent {\rm Let $(\bX, Y)$ be a random pair taking values in $\mathbb R^p \times \mathbb  R$. In the so-called single-index model, one has $Y=f^{\star}(\theta^{\star T}\bX)+\bW$, where $f^{\star}$ is an unknown univariate measurable function, $\theta^{\star}$ is an unknown vector in $\mathbb R^d$, and $W$ denotes a random noise satisfying $\mathbb E[\bW|\bX]=0$. The single-index model is known to offer a flexible way to model a variety of high-dimensional real-world phenomena. However, despite its relative simplicity, this dimension reduction scheme is faced with severe complications as soon as the underlying dimension becomes larger than the number of observations (``$p$ larger than $n$'' paradigm).  To circumvent this difficulty, we consider the single-index model estimation problem from a
sparsity perspective using a PAC-Bayesian approach. On the theoretical side, we offer a sharp oracle inequality, which is more powerful than the
best known oracle inequalities for other common procedures of single-index
recovery. The proposed method is implemented by means of the reversible jump Markov chain Monte Carlo technique and its performance is compared with that of standard procedures.
\medskip

\medskip
 
}

\end{abstract}

\setcounter{footnote}{0}
\section{Introduction}
Let $\mathcal D_n=\{(\bX_1,Y_1), \hdots, (\bX_n,Y_n)\}$ be a collection of
independent observations, distributed as a generic independent pair $(\bX, Y)$
taking values in $\mathbb R^p \times \mathbb  R$ and satisfying $\mathbb E Y^2 <
\infty$. Throughout, we let $\mathbf{P}$ be the distribution of $(\bX,
Y)$, so that the sample $\mathcal D_n$ is distributed according to $\mathbf
P^{\otimes n}$. In the regression function estimation problem, the goal is to
use
the data $\mathcal D_n$ in order to construct an estimate $r_n:\mathbb R^p \to
\mathbb R$ of the regression function $r(\bx)=\mathbb E[Y|\bX=\bx]$. In the
classical parametric linear model, one assumes
$$Y={\mathbf \theta}^{\star T}\bX+\bW,$$
where $\theta^{\star}=(\theta^{\star}_1, \hdots, \theta^{\star}_p)^T \in
\mathbb R^p$ and $\mathbb E[\bW|\bX]=0$. Here
$$r(\bx)=\theta^{\star T}\bx=\sum_{j=1}^p \theta^{\star}_j x_{j}$$
is a linear function of the components of $\bx=(x_1, \hdots, x_p)^T$. More
generally, we may define
\begin{equation}
\label{SI}
Y=f^{\star}(\theta^{\star T}\bX)+\bW,
\end{equation}
where $f^{\star}$ is an unknown univariate measurable function. This is the
celebrated
single-index model, which is recognized as a particularly useful variation of
the linear formulation and can easily be interpreted: The model changes only in
the direction $\theta^{\star}$, and the way it changes in this direction is
described by the function $f^{\star}$. This model has
applications to a variety of fields, such as discrete choice analysis in
econometrics and dose response models in biometrics, where high-dimensional
regression models are often employed. There are too many references to be
included here, but the monographs of McCullagh and Nelder \cite{MN} and Horowitz
\cite{Horo} together with the references \cite{Hardle,Ichimura,Delecroix,DJS,Lopez}
will
provide the reader with good introductions to the general subject area. 
\medskip

One of the main advantages of the single-index model is its supposed ability to
deal with the problem of high dimension (Bellman \cite{Bellman}). It is known that estimating the regression function is especially difficult whenever
the dimension $p$ of $\bX$ becomes large. As a matter of fact, the optimal mean
square convergence rate $n^{-2k/(2k+p)}$ for the estimation of a $k$-times
differentiable regression function converges to zero dramatically slowly if the
dimension $p$ is large compared to $k$. This leads to an unsatisfactory
accuracy of estimation for moderate sample sizes, and one possibility to
circumvent this problem is to impose additional assumptions on the regression
function. Thus, in particular, if $r(\bx)=f^{\star}(\theta^{\star T}\bx)$ holds
for every $\bx \in \mathbb R^p$, then the underlying structural dimension of the
model is 1 (instead of $p$) and the estimation of $r$ can hopefully be performed
easier.
In this regard, Ga{\"i}ffas and Lecu{\'e} show in \cite{GL} that the optimal
rate of convergence over the single-index model class is $n^{-2k/(2k+1)}$
(instead of $n^{-2k/(2k+p)}$), thereby answering a conjecture of Stone
\cite{Stone}.
\medskip

Nevertheless, practical estimation of the link function $f^{\star}$ and the
index $\theta^{\star}$ still requires a degree of statistical smoothing.
Perhaps the most common approach to reach this goal is to use a nonparametric
smoother (for instance, a kernel or a local polynomial method) to construct an
approximation $\hat f_n$ of $f^{\star}$, then substitute $\hat f_n$ into an
empirical version $R_n(\theta)$ of the mean square error $R(\theta)=\mathbb
E[Y-f(\theta^T\bX)]^2$, and finally choose $\hat \theta_n$ to minimize
$R_n(\theta)$ (see e.g.~H\"ardle, Hall and Ichimura \cite{Hardle} and Delecroix,
Hristache and
Patilea \cite{Delecroix} where the procedure is discussed in detail). The
rationale behind this type of two-stage approach, which is asymptotic in spirit, is that it produces a $\sqrt n$-consistent estimate of $\theta$, thereby devolving the difficulty to the
simpler problem of computing a good estimate for the one-dimensional function
$f^{\star}$. However, the relative simplicity of
this strategy is accompanied by severe difficulties (overfitting) when the dimension $p$ becomes larger than the number of observations $n$. Estimation in this setting (called ``$p$ larger than $n$'' paradigm) is generally acknowledged as an important challenge in contemporary statistics, see e.g.~the recent monograph of B\"uhlmann and van de Geer \cite{Sarah}. In fact, this drawback considerably reduces the ability of the single-index model to behave as an effective
dimension reduction technique.
\medskip

On the other hand, there is empirical evidence that many signals in
high-dimensional spaces admit a
sparse representation.  As an example, wavelet coefficients of images often
exhibit
exponential decay, and a relatively small subset of all wavelet coefficients
allow for a good approximation of the original image. Such signals have few
nonzero coefficients
and can therefore be described as sparse in the signal
domain (see for instance \cite{Bruckstein}). Similarly, recent advances in
high-throughput technologies---such as array comparative genomic hybridization---indicate that,
despite the huge dimensionality of problems, only a small number of genes may
play a role in determining the outcome and be required to create good predictors
(\cite{Veer} for instance). Sparse estimation is playing an increasingly
important role in the statistics and machine learning communities, and several
methods have recently been developed in both fields, which rely upon the notion
of sparsity (e.g.~penalty methods like the Lasso and Dantzig selector, see
\cite{Tibshirani,Candes, Bunea,Bickel} and the references therein). 
\medskip

In the present document, we consider the single-index model (\ref{SI}) from a
sparsity perspective, i.e., we assume that $\theta^{\star}$ has only a few
coordinates different from $0$. In the dimension reduction scenario we have in
mind, the ambient dimension $p$ can be very large, much larger than the sample
size $n$, but we believe that the representation is sparse, i.e., that very few
coordinates of $\theta^{\star}$ are nonzero. This assumption is helpful at
least for two reasons: If $p$ is large and the number of nonzero coordinates is
small enough, then the model is easier to interpret and its efficient estimation
becomes possible. Our setting is close in spirit of the approach of Cohen,  Daubechies, DeVore, Kerkyacharian and Picard \cite{Cohen}, who study approximation from queries of functions of the form $f(\theta^T \bx)$, where $\theta$ is approximately sparse (in the sense that it belongs to a weak-$\ell_{p}$ space). However, these authors dot not provide any statistical study of their model. Our modus operandi will rather rely on the so-called PAC-Bayesian
approach, originally developed in the classification context by Shawe-Taylor and
Williamson \cite{STW}, McAllester \cite{McAllester} and Catoni \cite{Catoni04,
Catoni07}. This strategy was further investigated for
regression by Audibert \cite{Audibert} and Alquier \cite{Alquier} and, more
recently, worked out in the sparsity framework by Dalalyan and Tsybakov
\cite{DT1,DT2} and Alquier and Lounici \cite{AL}. The main message of
\cite{DT1,DT2,AL} is that aggregation with a properly chosen prior is able to
deal nicely with the sparsity issue. Contrary to procedures such as the Lasso,
the
Dantzig selector and other penalized least square methods, which achieve fast rates under rather restrictive assumptions on the Gram matrix associated to
the predictors, PAC-Bayesian aggregation requires only minimal assumptions on the model.
Besides, it is computationally feasible even for a large $p$ and exhibits good
statistical performance.
\medskip

The paper is organized as follows. In Section \ref{sectiontheorique}, we first
set out some notation and introduce the single-index estimation procedure. Then we state our main result (Theorem
\ref{thm2}), which offers a sparsity oracle inequality more powerful than the
best known oracle inequalities for other common procedures of single-index
recovery. Section \ref{sectionexpe} is devoted to the practical implementation
of the estimate via a reversible jump Markov chain Monte Carlo (MCMC) algorithm, and to numerical experiments on both
simulated and real-life data sets. In order to preserve clarity, proofs have been postponed
to Section \ref{sectionproofs} and the description of the MCMC method in its
full
length is given in the Appendix Section \ref{sectionmcmc}.
\medskip

Note finally that our techniques extend to the case of multiple-index models, of
the form
$$Y=f^{\star}(\theta_1^{\star T}\bX, \hdots, \theta_m^{\star T}\bX)+\bW,$$
where the underlying structural dimension $m$ is supposed to be larger than 1
but substantially smaller than $p$. However, to keep things simple, we let $m=1$
and leave the reader the opportunity to adapt the results to the more general
situation $m\geq 1$.
\section{Sparse single-index estimation}
\label{sectiontheorique}
We start this section with some notation and basic requirements.
\subsection{Notation}
Throughout the document, we suppose that the recorded data $\mathcal D_n$ is
generated according to the single-index model (\ref{SI}). More precisely, for
each $i=1, \hdots, n$,  
\begin{equation*}
\label{SI2}
Y_i=f^{\star}(\theta^{\star T}\bX_i)+\bW_i,
\end{equation*}
where $f^\star$ is a univariate measurable function, $\theta^{\star}$ is a
$p$-variate
vector, and $\bW_1, \hdots, \bW_n$ are independent copies of $W$. We emphasize that it is implicitly assumed that the observations are drawn according to the true model under study. Therefore, this casts our problem in a nonparametric setting rather than in a classical PAC-Bayesian one. 
\medskip

Recall that, in model (\ref{SI}), $\mathbb E[\bW|\bX]=0$ and, consequently, that $\mathbb
E \bW=0$. However, the distribution of $\bW$ (in particular, the variance)
may depend on $\bX$. We shall not precisely specify this dependence, and will
rather require the following condition on the distribution of $\bW$. 
\medskip

\noindent\textbf{Assumption $\mathbf N$.}  There exist two positive constants
$\sigma$ and $L$ such that, for
all integers $k \geq 2$, 
$$ \mathbb{E}\left[|\bW|^{k}\,|\,\bX\right] \leq
\frac{k!}{2}\sigma^{2}L^{k-2}.$$
\medskip

Observe that Assumption $\mathbf N$ holds in particular if $\bW=\Phi(\bX)\varepsilon$,
where $\varepsilon$ is a standard Gaussian random variable independent of $\bX$ and
$\Phi(\bX)$ is almost surely bounded. 
\medskip

Let $\|\theta\|_1$ denote the $\ell_1$-norm of the vector $\theta=(\theta_1,
\hdots, \theta_p)^T$, i.e., $\|\theta\|_1=\sum_{j=1}^p|\theta_j|$. Without loss of generality, it will be assumed throughout the document that the index $\theta^{\star}$ belongs to $\mathcal{S}^{p}_{1,+}$,
where $\mathcal{S}^{p}_{1,+}$ is the set of all $\theta\in\mathbb{R}^{p}$ with
$\|\theta\|_{1}=1$ and $\theta_{j(\theta)}>0$, where $j(\theta)$ is the
smallest $j\in\{1,\hdots,p\}$ such that $\theta_{j}\neq 0$.
\medskip

We will also require that the random variable $\bX$ is
almost surely bounded by a constant which, without loss of generality, can be
taken equal to 1.  Moreover, it will also be assumed that the link function
$f^{\star}$ is bounded by some known positive constant $C$. Thus, denoting by
$\|\bX\|_{\infty}$ the supremum norm of $\bX$ and by $\|f^{\star}\|_{\infty}$
the functional
supremum norm of $f^{\star}$ over $[-1,1]$, we set:
\medskip

\noindent\textbf{Assumption $\mathbf B$.} The condition $\|\bX\|_{\infty} \leq 1$ holds almost surely and there exists a positive constant $C$ larger than 1 such that $\|f^{\star}\|_{\infty} \leq C.$
\begin{rem}
To keep a sufficient degree of clarity, no attempt was made to optimize the constants. In particular, the requirement $C\geq1$ is purely technical. It can easily be removed by replacing $C$ by $\max[C,1]$ throughout the document.
\end{rem}
In order to approximate the link function $f^{\star}$, we shall use the vector
space $\mathcal F$ spanned by a given countable dictionary of measurable
functions $\{\varphi_{j}\}_{j=1}^{\infty}$. Put differently, the approximation
space $\mathcal F$ is the set of (finite) linear combinations of functions of
the dictionary. Each $\varphi_j$ of the collection is assumed
to be defined on $[-1,1]$ and to take values in $[-1,1]$. To avoid getting into too
much technicalities, we will also assume that each $\varphi_j$ is differentiable and such that, for some positive constant $\ell$, $\|\varphi_j'\|_\infty \leq
\ell j$. This assumption is satisfied by the (non-normalized) trigonometric system
$$\varphi_1(t)=1,\, \varphi_{2j}(t)=\cos(\pi j t),\,
\varphi_{2j+1}(t)= \sin (\pi j t),\quad j=1, 2, \hdots$$

Finally, for any measurable $f:\mathbb{R}^{p}\to\mathbb {R}$ and
$\theta\in\mathcal{S}^{p}_{1,+}$, we let
$$R(\theta,f)=\mathbb{E} \left[ \left(Y- f(\theta^{T}\bX)\right)^{2}\right]$$
and denote by
$$R_{n}(\theta,f)=\frac{1}{n}\sum_{i=1}^{n} \left(Y_{i}- f(\theta^T\bX_{i})
\right)^{2}$$
the empirical counterpart of $R(\theta, f)$ based on the sample $\mathcal D_n$.
\subsection{Estimation procedure}
\label{PLUS}
We are now in a position to describe our estimation procedure. The method which
is presented here is inspired by the approach developed by Catoni in
\cite{Catoni04,Catoni07}. It strongly relies on the choice of a
probability measure $\pi$ on $\mathcal S_{1,+}^p\times \mathcal F$, called the
prior, which in our framework should enforce the sparsity properties of the
target regression function. With this objective in mind, we first let
$$ \mbox{d}\pi(\theta,f) = \mbox{d}\mu(\theta) \mbox{d}\nu(f),$$
i.e., we assume that the distribution over the indexes is independent
of the distribution over the link functions.
With respect to the parameter $\theta$, we put
\begin{equation}
\label{1plus}
\mbox{d} \mu(\theta) = \frac{\displaystyle
\sum_{i=1}^{p}10^{-i}\sum_{I\subset\{1,\hdots,p\},|I|=i} {p\choose i}^{-1}
                           \mbox{d}\mu_{I}(\theta)}{1-(\frac{1}{10})^{p}},
\end{equation}
where $|I|$ denotes the cardinality of $I$ and $\mbox{d}\mu_{I}(\theta)$ is the
uniform probability measure on the set
$$\mathcal{S}^{p}_{1,+}(I)=\{\theta=(\theta_1, \hdots, \theta_p)\in\mathcal{S}^{p}_{1,+}: \theta_{j}=0
\mbox{ if and only if } j\notin I\}.$$
We see that $\mathcal{S}^{p}_{1,+}(I)$ may be interpreted as the set of ``active''
coordinates in the single-index regression of $Y$ on $\bX$, and note that the prior on
$\mathcal S_{1,+}^p$ is a convex combination of uniform probability measures on
the subsets $\mathcal S_{1,+}^p(I)$. The weights of this combination depend only
on the size of the active coordinate subset $I$. As such, the value $|I|$ characterizes the
sparsity of the model: The smaller $|I|$, the smaller the number of variables involved in the model. The factor $10^{-i}$ penalizes models of high dimension, in accordance with the sparsity idea.
\medskip

The choice of the prior $\nu$ on $\mathcal F$ is more involved. To
begin with, we define, for any positive integer $M\leq n$ and all $\Lambda>0$, 
\begin{equation*}
\mathcal B_{M}(\Lambda) =\left \{(\beta_{1},\hdots,\beta_{M})\in\mathbb{R}^{M}:
\sum_{j=1}^M j|\beta_{j}| \leq \Lambda \text{ and } \beta_{M} \neq 0 \right\}.
\end{equation*}
Next, we let $\mathcal F_M(\Lambda)\subset \mathcal F$
be the image of $\mathcal B_M(\Lambda)$ by the map 
\begin{equation*}
\begin{array}{clcl}
{\Phi}_{M}& :\mathbb R^M &\to & \mathcal F\\
 & (\beta_{1},\hdots,\beta_{M}) &\mapsto &\sum_{j=1}^{M} \beta_{j}\varphi_{j}.
\end{array}
\end{equation*}
It is worth pointing out that, roughly, Sobolev spaces are well approximated by $\mathcal F_M(\Lambda)$  as $M$ grows (more on this in Subsection \ref{zzz}). Finally, we define $\nu_{M}(\mbox{d}f)$ on the set $\mathcal
F_M(C+1)$ as the image of the
uniform measure on $\mathcal B_M(C+1)$ induced by the map ${\Phi}_{M}$, and
take
\begin{equation}
\label{2plus}
\mbox{d} \nu(f) =  \frac{\displaystyle \sum_{M=1}^{n}10^{-M}
                           \mbox{d}\nu_{M}(f)}{1-(\frac{1}{10})^{n}} .
\end{equation}
Some comments are in order here. First, we note that the prior $\pi$ is defined on $\mathcal S_{1,+}^p\times \mathcal F_n(C+1)$ endowed with its canonical Borel $\sigma$-field. The choice of $C+1$ instead of $C$ in the definition of the prior support is essentially technical. This bound ensures that when the target $f^{\star}$ belongs to $\mathcal F_n(C)$, then a small ball around it is contained in $\mathcal F_n(C+1)$. It could be safely replaced by $C+u_n$, where $\{u_n\}_{n=1}^{\infty}$ is any positive sequence vanishing sufficiently slowly as $n\to \infty$. Next, the integer $M$ should be interpreted as a measure of the
``dimension'' of the function $f$---the larger $M$, the more complex the
function---and the prior $\nu$ adapts again to the sparsity idea
by penalizing large-dimensional functions $f$. The coefficients $10^{-i}$ and $10^{-M}$ which appear in (\ref{1plus}) and (\ref{2plus}) show that more complex models have a geometrically decreasing influence. Note however that the value 10, which has been chosen because of its good practical results, is somehow arbitrary. It could be, in all generality, replaced by a more general coefficient $\alpha$ at the price of a more technical analysis. Finally, we observe that, for
each $f = \sum_{j=1}^M \beta_j \varphi_j \in \mathcal F_M(C+1)$,
\begin{align*}
\|f\|_{\infty} \leq \sum_{j=1}^M |\beta_j| \leq C+1.
\end{align*}
Now, let $\lambda$ be a positive real number, called the inverse temperature parameter
hereafter. The estimates $\hat{\theta}_{\lambda}$ and $\hat{f}_{\lambda}$ of
$\theta^{\star}$ and $f^{\star}$, respectively, are simply obtained by randomly drawing
$$ (\hat{\theta}_{\lambda},\hat{f}_{\lambda}) \sim \hat{\rho}_{\lambda},$$
where $\hat{\rho}_{\lambda}$ is the so-called Gibbs posterior distribution over $\mathcal S_{1,+}^p\times \mathcal F_n(C+1)$, defined by the probability density 
$$ \frac{\mbox{d}\hat{\rho}_{\lambda}}{\mbox{d}\pi}(\theta, f)= \frac{\exp\left
[-\lambda R_{n}(\theta,f)\right]}{\displaystyle \int \exp\left [-\lambda
R_{n}(\theta,f)\right]\mbox{d}\pi(\theta,f)} .$$
[The notation ${\mbox{d}\hat{\rho}_{\lambda}}/{\mbox{d}\pi}$ means the density
of $\hat{\rho}_{\lambda}$ with respect to $\pi$.]
The estimate $(\hat{\theta}_{\lambda},\hat{f}_{\lambda})$ has a simple
interpretation. Firstly, the level of significance of each pair $(\theta, f)$
is assessed via its least square error performance on the data $\mathcal D_n$.
Secondly, a Gibbs distribution with respect to the prior $\pi$ enforcing
those pairs $(\theta, f)$ with the most empirical significance is assigned on
the space $\mathcal S_{1,+}^p \times \mathcal F_n(C+1)$. Finally, the resulting
estimate is
just a random realization (conditional to the data) of this Gibbs  posterior
distribution.
\subsection{Sparsity oracle inequality}
\label{zzz}
For any
$I\subset\{1,\hdots,p\}$ and any positive integer $M\leq n$, we set
$$ \left(\theta^{\star}_{I,M},f^{\star}_{I,M}\right) \in
\arg\min_{(\theta,f)\in\mathcal{S}^{p}_{1,+}(I)\times\mathcal{F}_{M}(C)}
               R(\theta,f). $$
At this stage, it is very important to note that, for each $M$, the infimum $f^{\star}_{I,M}$ is defined on $\mathcal F_{M}(C)$, whereas the prior charges a slightly bigger set, namely $\mathcal F_{M}(C+1)$.             
\medskip             
               
The main result of the paper is the following theorem. Here and everywhere, the
wording ``with probability $1-\delta$'' means the probability evaluated
with respect to the distribution $\mathbf P^{\otimes n}$ of the data $\mathcal
D_n$ {\it and} the conditional probability measure $\hat \rho_{\lambda}$. Recall that $\ell$ is a positive constant such that $\|\varphi_j'\|_\infty \leq
\ell j$. 
\begin{theo}
\label{thm2}
Assume that Assumption $\mathbf N$ and Assumption $\mathbf B$ hold.  Set 
$w=8(2C+1)\max[L,2C+1]$ and take 
\begin{equation}
\label{equationlambda}
\lambda = \frac{n}{w+2\left[(2C+1)^{2}+4\sigma^{2}\right]}.
\end{equation}
Then, for all $\delta\in\,]0,1[$, with probability at least
$1-\delta$
we have
\begin{align*}
R(\hat{\theta}_{\lambda},\hat{f}_{\lambda})- R(\theta^{\star},f^{\star}) &\leq
\Xi\inf_{
\tiny{\begin{array}{c} I\subset\{1,\hdots,p\}
                                 \\ 1\leq M \leq n \end{array}}} \Bigg\{
           R(\theta^{\star}_{I,M},f^{\star}_{I,M}) -
R(\theta^{\star},f^{\star}) \\
&\qquad \qquad  + \frac{ M\log(Cn) + |I|\log(pn)+
\log\left(\frac{2}{\delta}\right)}{n}\Bigg\},
\end{align*}
where $\Xi$ is a positive constant, function of $L$, $C$, $\sigma$ and $\ell$ only.
\end{theo}
\begin{rem}
Interestingly enough, analysis of the estimate $(\hat{\theta}_{\lambda},\hat{f}_{\lambda})$ is still possible when Assumption $\mathbf N$ is not satisfied. Indeed, even if Bernstein's inequality (see Lemma \ref{lemmemassart}) is not valid, a recent paper by Seldin, Cesa-Bianchi, Laviolette, Auer, Shawe-Taylor and Peters \cite{Seldin} provides us with a nice alternative inequality assuming less restrictive assumptions. However, we would then suffer a loss in the upper bound of Theorem \ref{thm2}. It is also interesting to note that recent results by Audibert and Catoni \cite{AudibertCatoni} allow the study of PAC-Bayesian estimates without Assumption {\bf N}. However, the results of these authors are valid for linear models only, and it is therefore not clear to what extent their technique can be transposed to our setting. 
\end{rem}

Theorem \ref{thm2} can be given a simple interpretation. Indeed, we see that if there is a
``small'' $I$ and a ``small'' $M$ such that
$R(\theta^{\star}_{I,M},f^{\star}_{I,M})$ is close to
$R(\theta^{\star},f^{\star})$, then
 $R(\hat{\theta}_{\lambda},\hat{f}_{\lambda})$ is also close
to $R(\theta^{\star},f^{\star})$ up to terms of order $1/n$. However, if no
such $I$ or $M$ exists, then one of the terms $M\log(Cn)/n$ and
$|I|\log(pn)/n$ starts to dominate, thereby deteriorating the general
quality of the bound. A good approximation with a ``small'' $I$ is typically
possible
when $\theta^{\star}$ is sparse or, at least, when it can be approximated by a
sparse
parameter. On the other hand, a good approximation with a ``small'' $M$ is
possible if $f^{\star}$ has a sufficient degree of regularity.
\medskip

To illustrate the latter remark,  assume for instance that
$\{\varphi_j\}_{j=1}^{\infty}$ is the (non-normalized) trigonometric system and
suppose that the target $f^{\star}$ be\-longs to the So\-bo\-lev ellipsoid, defined by
$$\mathcal W\left(k,\frac{{6}C^2}{\pi^2}\right)= \left\{ f\in L_{2}([-1,1]):
f=\sum_{j=1}^{\infty}\beta_j \varphi_j \mbox{ and } 
     \sum_{j=1}^{\infty} j ^{2k}
\beta_{j}^{2}  \leq \frac{{6}C^2}{\pi^2} \right\}
$$
for some unknown regularity parameter $k\geq 2$ (see, e.g., Tsybakov
\cite{Tsybakov}). Observe that, in this context, the approximation sets
$\mathcal F_M(C+1)$ take the form
\begin{align*}
& \mathcal F_M(C+1)\\
& \quad = \left\{ f\in L_{2}([-1,1]): f=\sum_{j=1}^{M}\beta_j
\varphi_j, \sum_{j=1}^{M} j|\beta_{j}|  \leq C+1 \mbox{ and } \beta_M\neq 0\right\}.
\end{align*}
It is important to note that the regularity parameter $k$ is assumed to be unknown,
and this casts our results in the so-called adaptive setting. The following additional assumption will be needed:
\medskip

\noindent\textbf{Assumption $\mathbf D$.} The random variable $\theta^{\star T}\bX$ has a probability density on
$[-1,1]$,
bounded from above by a positive constant $B$.
\medskip

Last, we let
$I^{\star}$ be the set $I$ such that
$\theta^{\star}\in\mathcal{S}_{1,+}^{p}(I)$ and set
$\|\theta^\star\|_{0}=|I^\star|$.
\begin{cor}
\label{thm3}
Assume that Assumption $\mathbf N$, Assumption $\mathbf B$ and Assumption $\mathbf D$ hold. Suppose also that $f^{\star}$ be\-lon\-gs to the So\-bo\-lev
ellipsoid $ \mathcal W(k,{6}C^2/\pi^2)$, where the real number $k\geq 2$ is an (unknown)
regularity parameter. Set  $w=8(2C+1)\max[L,2C+1]$ and take 
$\lambda$ as in (\ref{equationlambda}).
Then, for all $\delta\in\,]0,1[$, with probability at least
$1-\delta$
we have
\begin{align}
& R(\hat{\theta}_{\lambda},\hat{f}_{\lambda}) -
R(\theta^{\star},f^{\star})\nonumber
\\
&\leq \Xi' \left\{
\left(\frac{\log(Cn)}{n}\right)^{\frac{2k}{2k+1}} +
\frac{\|\theta^\star\|_{0} \log(pn)}{n}+\frac{
\log\left(\frac{2}{\delta}\right)}{n}\right\},\label{kx}
\end{align}
where $\Xi'$ is a positive constant, function of $L$, $C$, $\sigma$,
$\ell$ and $B$ only.
\end{cor}

As far as we are aware, all existing methods achieving rates of convergence similar to the ones provided by Corollary \ref{thm3} are valid in an asymptotic setting only ($p$ fixed and $n \to \infty$).
The strength of Corollary \ref{thm3} is to provide a finite sample bound and to show that our estimate still behaves well in a nonasymptotic situation if the intrinsic dimension (i.e., the sparsity) is small with respect to $n$. To understand this remark, just assume that $p$ is a function of $n$ such that $p\to \infty$ as $n \to \infty$. Whereas a classical asymptotic approach cannot say anything clever about this situation, our bounds still provide some information, provided the model is sparse enough (i.e., $\|\theta^{\star}\|_{0}$ is sufficiently small with respect to $n$).
\medskip

We see that, asymptotically ($p$ fixed and $n \to \infty$), the leading term on the right-hand side of
inequality (\ref{kx}) is
$(\log(n)/n)^{\frac{2k}{2k+1}}$.
This is the minimax rate of convergence over a Sobolev class, up to a $\log(n)$ factor. However,
when $n$ is ``small'' and $\theta^{\star}$ is not sparse (i.e.,
$\|\theta^\star\|_{0}$ is not
``small''), the term $\|\theta^\star\|_{0}\log(pn)/n$
starts to emerge and cannot be neglected. Put differently, in large dimension,
the estimation of $\theta^ {\star}$
itself is a problem---this phenomenon is not taken into account by asymptotic
studies.
\medskip

It is worth mentioning that the approach developed in the present article does not offer any guarantee on the point of view of variable (feature) selection. To reach this objective, an interesting route to follow is  the sufficient dimension reduction (SDR) method proposed by Chen, Zou and Cook \cite{Chen}, which can be applied to the 
single-index model to estimate consistently the parameter $\theta^\star$ and perform variable selection in a sparsity framework. Note however that such results require strong assumptions on the distribution of the data.
\medskip

Finally, it should be stressed that the choice of $\lambda$ in Theorem
\ref{thm2} and Corollary \ref{thm3} is not the best possible and may eventually
be improved, at the price of a more technical analysis however.

\section{Implementation and numerical results}
\label{sectionexpe}
A series of experiments was conducted, both on simulated and real-life data sets, in order to assess the practical capabilities of the proposed method and compare its performance with that of standard procedures. Prior to analysis, we first need to discuss its concrete implementation, which has been carried out via  a Markov Chain Monte Carlo (MCMC) method. 
\subsection{Implementation via reversible jump MCMC}
The use of MCMC methods has become a popular way to compute Bayesian estimates.  For an introduction to the domain, one should refer to the comprehensive monograph of Marin and Robert \cite{Christian2} and the references
therein. Importantly, in this computational framework,  an adaptation of the well-known Hastings-Metropolis algorithm to the case where the posterior distribution gives mass to several
models of different dimensions was proposed by Green \cite{RJMCMC} under the
name Reversible Jump MCMC (RJMCMC) method.  In the PAC-Bayesian setting,  MCMC procedures were first 
considered by Catoni \cite{Catoni04}, whereas Dalalyan
and Tsybakov \cite{DT1,DT2} and Alquier and Lounici \cite{AL} explore their practical implementation in the sparse context using Langevin Monte Carlo and RJMCMC, respectively. Regarding the single-index model,  MCMC algorithms were used to compute Bayesian estimates by Antoniadis, Gr\'egoire and
McKeague \cite{Antoniadis} and, more recently, by Wang \cite{Wang}, who develop a fully Bayesian method to analyse the single-index model. Our implementation technique is close in spirit to the one of Wang \cite{Wang}.
\medskip

As a starting point for the approximate computation of our estimate, we used the RJMCMC method of Green \cite{RJMCMC}, which is in fact an adaptation of the
Hastings-Metropolis algorithm to the case where the objective posterior probability
distribution (here, $\hat{\rho}_{\lambda}$) assigns mass to several different
models. The idea is to start from an initial given pair $(\theta^{(0)},f^{(0)}) \in \mathcal S_{1,+}^p \times \mathcal F_n(C+1)$ and then, at each step, to iteratively compute
$(\theta^{(t+1)},f^{(t+1)})$ from $(\theta^{(t)},f^{(t)})$ via the following chain of rules:
\begin{itemize}
\item Sample a random pair $(\tau^{(t)},h^{(t)})$ according to some proposal conditional density
$k_{t}(\,.\,|(\theta^{(t)},f^{(t)}))$ with respect to the prior $\pi$;
\item Take
$$ (\theta^{(t+1)},f^{(t+1)}) = \left\{
\begin{array}{l l}
(\tau^{(t)},h^{(t)}) & \text{ with probability } \alpha_{t}
\\
(\theta^{(t)},f^{(t)}) & \text{ with probability }
1-\alpha_{t},
\end{array}
\right. $$
where 
$$\alpha_{t}= \min\left(1,\frac{ \frac{\mbox{\footnotesize d}\hat{\rho}_{\lambda}}
                    {\mbox{\footnotesize d}\pi}(\tau^{(t)},h^{(t)})\times
          k_t\left ((\theta^{(t)},f^{(t)})|(\tau^{(t)},h^{(t)})\right)}
          {\frac{\mbox{\footnotesize d}\hat{\rho}_{\lambda}}{\mbox{\footnotesize d}\pi}(\theta^{(t)},f^{(t)})\times
                    k_t\left((\tau^{(t)},h^{(t)})|(\theta^{(t)},f^{(t)})\right)}\right).$$
\end{itemize}
This protocol ensures that the sequence $\{(\theta^{(t)},f^{(t)})\}_{t=0}^{\infty}$ is a Markov chain with invariant
probability distribution $\hat{\rho}_{\lambda}$ (see e.g.~Marin and Robert \cite{Christian2}). A usual choice is to take $k_{t}\equiv k$, so that the Markov chain is homogeneous. However, in our context, it is more convenient to let $k_{t}=k_{1}$ if $t$ is odd and $k_{t}=k_{2}$ if $t$ is even. Roughly, the effect of $k_{1}$ is to modify the index $\theta^{(t)}$ while $k_{2}$ will essentially act on the link function $f^{(t)}$. While the ideas underlying the proposal densities $k_{1}$ and $k_{2}$ are quite simple, a precise description
in its full length turns out to be more technical. Thus, in order to preserve the readability of
the
paper, the explicit construction of $k_1$ and $k_2$ has been postponed to the Appendix Section \ref{sectionmcmc}.
\medskip

From a theoretical point of view, it is clear that the implementation of our method requires knowledge of the constant $C$ (the upper bound on $\|f^\star\|_{\infty}$). A too small $C$ will result in a smaller model, which is unable to perform a good approximation. On the other hand, a larger $C$ induces a poor bound in Theorem 2.1. In practice, however, the influence of $C$ turns out to be secondary compared to the impact of the parameter $\lambda$. Indeed, it was found empirically that a very large choice of $C$ (e.g., $C=10^{100}$) does not deteriorate the overall quality of the results, as soon as $\lambda$ is  appropriately chosen. This is the approach that was followed in the experimental testing process.
\medskip

Besides, the time for the Markov chains to converge depends strongly on the ambient dimension $p$ and the starting point of the simulations. When the dimension is small (typically, $p\leq 10$), the chains converge fast and any value may be chosen as a starting point. In this case, we let the MCMC run $1000$ steps and obtained satisfying results. On the other hand, when the dimension is larger (typically, $p>10$), the convergence is very slow, in the sense that $R_n (\theta^{(t)},f^{(t)})$ takes a very long time to stabilize. However, using as a starting point for the chains the preliminary estimate $\hat{\theta}_{\mbox{\footnotesize HHI}}$ (see below) significantly reduces the number of steps needed to reach  convergence---we let the chains run $5000$ steps in this context. Nevertheless, as a general rule, we encourage the users to inspect the convergence of the chains by checking if $R_n (\theta^{(t)},f^{(t)})$ is stabilized, and to run several chains starting from different points to avoid their attraction into local minima.
\subsection{Simulation study}
In this subsection, we illustrate the finite sample performance of the presented estimation method on  three synthetic data sets and compare its predictive capabilities with those of three standard statistical procedures. In all our experiments, we took as dictionary the (non-normalized) trigonometric system $\{\varphi_j\}_{j=1}^{\infty}$ and denote accordingly the resulting regression function estimate defined in Section \ref{sectiontheorique} by $\hat{F}_{\mbox{\footnotesize Fourier}}$. In accordance with the order of magnitude indicated by the theoretical results, we set $\lambda=4n$. This choice can undoubtedly be improved a bit but, as the numerical results show, it seems sufficient for our procedure to be fairly competitive.
\medskip

The tested competing methods are the Lasso (Tibshirani
\cite{Tibshirani}), the standard regression kernel estimate (Nadaraya \cite{Nadaraya1,Nadaraya2} and Watson \cite{Watson}, see also Tsybakov \cite{Tsybakov}), and the estimation strategy discussed in H\"ardle, Hall and Ichimura \cite{Hardle}. While the procedure of H\"ardle, Hall and Ichimura is specifically tailored for single-index models, the Lasso is designed to deal with the estimation of sparse linear
models. On the other hand, the nonparametric kernel method is one of the best options when no obvious assumption (such as the single-index one) can be made on the shape of the targeted regression function. 
\medskip

We briefly recall that, for a linear model of the form $Y=\theta^{\star T}\bX+\bW$, the Lasso estimate takes the form $\hat{F}_{\mbox{\footnotesize Lasso}}(\bx)= \hat{\theta}_{\mbox{\footnotesize Lasso}}^{T}\bx$, where
$$ \hat{\theta}_{\mbox{\footnotesize Lasso}} \in \arg\min_{\theta \in\mathbb{R}^{p}} \left\{
\frac{1}{n}\sum_{i=1}^{n} \left(Y_i - \theta^{T} \mathbf{X}_{i}\right)^{2}
+ \xi\sum_{j=1}^{p} |\theta_{j}|\right\}$$
and $\xi>0$ is a regularization parameter. Theoretical results (see e.g.~Bunea, Tsybakov and Wegkamp \cite{Bunea}) indicate that $\xi$ should be of the order
$\xi^{\star}=\sigma\sqrt{\log(p)/n}$. Throughout, $\sigma$ is assumed to be known, and we let $\xi=\xi^{\star}/3$, since this choice is known to give good practical results. The Nadaraya-Watson kernel estimate will be denoted by $\hat{F}_{\mbox{\footnotesize NW}}$. It is defined by
$$ \hat{F}_{\mbox{\footnotesize NW}}(\bx) = \frac{\sum_{i=1}^{n} Y_{i} K_h(\bx-\mathbf{X}_{i})}
{\sum_{i=1}^{n} K_h(\bx-\mathbf{X}_{i})}
$$
for some nonnegative kernel $K$ on $\mathbb{R}^{p}$ and $K_{h}(\mathbf z)=K(\mathbf z/h)/h$. In the experiments, we let
$K$ be the Gaussian kernel $K(\mathbf z)=\exp(-\mathbf z^T\mathbf z)$ and chose the smoothing parameter $h$ via a classical leave-one-out procedure on the grid
$\mathcal G=\{0.75^{k},k=0,\hdots,\lfloor\log(n)\rfloor\}$, see, e.g., Gy\"orfi, Kohler, Krzy\.zak and Walk \cite{Laciregressionbook} (notation $\lfloor.\rfloor$ stands for the floor function). Finally, the estimation procedure advocated in H\"ardle, Hall and Ichimura \cite{Hardle} takes the form
$$ \hat{F}_{\mbox{\footnotesize HHI}}(\bx) = \frac
{\sum_{i=1}^{n} Y_i G_{\hat{h}}\left(\hat{\theta}^{T}_{\mbox{\footnotesize HHI}}(\bx-\bX_{i})\right)} {\sum_{i=1}^{n} G_{\hat{h}}\left(\hat{\theta}^{T}_{\mbox{\footnotesize HHI}}(\bx-\bX_{i})\right)}
 $$
for some kernel $G$ on $\mathbb{R}$, with $G_{h}(\mathbf z)=G(\mathbf z/h)/h$ and
$$ \left(\hat{h},\hat{\theta}_{\mbox{\footnotesize HHI}}\right) 
    \in \arg\min_{h>0,\theta\in \mathbb R^p} \sum_{i=1}^{n}
          \left[Y_{i} - \frac
{\sum_{j\neq i} Y_j G_{h}\left(\theta^{T}(\bX_{j}-\bX_{i})\right)}
{\sum_{j\neq i} G_{h}\left(\theta^{T}(\bX_{j}-\bX_{i})\right)}\right]^{2}.$$
All calculations were performed with the Gaussian kernel. We used the grid $\mathcal G$ for the optimization with respect to $h$, whereas the best search for $\theta$ was implemented via a pathwise coordinate optimization. 
\medskip

The various methods were tested for the general regression model 
$$ Y_i= F(\mathbf{X}_i) + {W}_i,\quad  i=1, \hdots, n,$$
for three different choices of $F$ (single-index or
not) and two values of $n$, namely $n=50$ and $n=100$. In each of these models, the observations $\mathbf{X}_i$ take values in $\mathbb{R}^{p}$, with $p=10$ and $p=50$, and have independent components uniformly distributed on $[-1,1]$. The noise variables ${W}_1, \hdots, W_n$ are independently distributed according to a Gaussian $\mathcal{N}(0,\sigma^{2})$, with $\sigma=0.2$. It is worth pointing out that for $n=50$ and $p=50$, $p$ and $n$ are
of the same order, which means that the setting is nonasymptotic. It is essentially in this case that the use of estimates tailored to sparsity, which reduce the variance, is expected to improve the performance over generalist methods. On the other hand, the situation $n=100$ and $p=10$ is less difficult and mimics the asymptotic setting.
\medskip

The three examined functions $F(\bx)$, for $\bx=(x_1,\hdots,x_{p})$, were the following ones:
\begin{itemize}
 \item[][{\bf Model 1}] A linear model $F_{\mbox{\footnotesize Linear}}(\bx)=2 \theta^{\star {T}}\bx$.
 \item[][{\bf Model 2}] A single-index function $F_{\mbox{\footnotesize SI}}(\bx) = 2(\theta^{\star T}\bx)^{2} + \theta^{\star T}\bx$.
 \item[][{\bf Model 3}] A purely nonparametric model
           $F_{\mbox{\footnotesize NP}}(\bx) = 2|x_{2}|\sqrt{|x_{1}|} -x_{3}^{3}$,
\end{itemize}
where, in the first and second model, $\theta^{\star}=(0.5,0.5,0,\hdots,0)^{T}$. Thus, in [{\bf Model 1}] and [{\bf Model 2}], even if the ambient dimension is large, the intrinsic dimension of the model is in fact equal to 2. 
\medskip

For each experiment, a learning set of size $n$ was generated to compute the estimates and their performance, in terms of mean square prevision error, was evaluated on a separate test set of the same size. The results are shown in Table \ref{table1} ($p=10)$ and Table \ref{table1bis} ($p=50)$. As each experiment was repeated $20$ times, these tables report the median, the mean and the standard deviation (s.d.)~of the prevision error of each procedure.
\begin{table}[t!]
\begin{center}
\begin{tabular}{|p{1.6cm}|p{1.4cm}||p{1.7cm}|p{1.7cm}|p{1.7cm}|p{1.7cm}|p{1.3cm}
|}
\hline
\hline
 $n=50$  & $p=10$       & $\hat{F}_{\mbox{\footnotesize Fourier}}$ & $\hat{F}_{\mbox{\footnotesize HHI}}$
                             & $\hat{F}_{\mbox{\footnotesize Lasso}}$ & $\hat{F}_{\mbox{\footnotesize NW}}$ \\
\hline \hline
$F_{\mbox{\footnotesize Linear}}$    & median & 0.061 & 0.063 & {\bf 0.046} & 0.293 \\
                & mean   & 0.061 & 0.063 & {\bf 0.047} & 0.290 \\
                & s.d.   & 0.016 & 0.014 & 0.011 & 0.063 \\
\hline
$F_{\mbox{\footnotesize SI}}$    & median & {\bf 0.050} & 0.067 & 0.307 & 0.198 \\
                & mean   & {\bf 0.069} & 0.080 & 0.338 & 0.208 \\
                & s.d.   & 0.081 & 0.057 & 0.082 & 0.072 \\
\hline
$F_{\mbox{\footnotesize NP}}$    & median & 0.375 & 0.405 & 0.830 & {\bf 0.354} \\
                & mean   & 0.402 & 0.407 & 0.890 & {\bf 0.336} \\
                & s.d.   & 0.166 & 0.110 & 0.176 & 0.006 \\
\hline
\end{tabular}
\vspace*{-6pt}
\begin{tabular}{|p{1.6cm}|p{1.4cm}||p{1.7cm}|p{1.7cm}|p{1.7cm}|p{1.7cm}|p{1.3cm}
|}
\hline $n=100$  &  $p=10$       & $\hat{F}_{\mbox{\footnotesize Fourier}}$ & $\hat{F}_{\mbox{\footnotesize HHI}}$
                             & $\hat{F}_{\mbox{\footnotesize Lasso}}$ & $\hat{F}_{\mbox{\footnotesize NW}}$ \\
\hline \hline
$F_{\mbox{\footnotesize Linear}}$    & median & 0.053 & 0.051 & {\bf 0.042} & 0.227 \\
                & mean   & 0.056 & 0.050 & {\bf 0.043} & 0.237 \\
                & s.d.   & 0.011 & 0.006 & 0.004 & 0.044 \\
\hline
$F_{\mbox{\footnotesize SI}}$    & median & {\bf 0.047} & 0.052 & 0.332 & 0.209 \\
                & mean   & {\bf 0.049} & 0.053 & 0.337 & 0.218 \\
                & s.d.   & 0.009 & 0.012 & 0.063 & 0.045 \\
\hline
$F_{\mbox{\footnotesize NP}}$    & median & {\bf 0.305} & 0.343 & 0.793 & 0.333 \\
                & mean   & {\bf 0.321} & 0.338 & 0.833 & 0.324 \\
                & s.d.   & 0.092 & 0.042 & 0.145 & 0.041 \\
\hline
\hline
\end{tabular}
\caption{\textsf{Numerical results for the simulated data, with $n=50$ and $n=100$,
$p=10$.
The characters in bold indicate the best performance.}}
\label{table1}
\end{center}
\end{table}
\begin{table}[t!]
\begin{center}
\begin{tabular}{|p{1.6cm}|p{1.4cm}||p{1.7cm}|p{1.7cm}|p{1.7cm}|p{1.7cm}|p{1.3cm}
|}
\hline
\hline
 $n=50$  & $p=50$       & $\hat{F}_{\mbox{\footnotesize Fourier}}$ & $\hat{F}_{\mbox{\footnotesize HHI}}$
                             & $\hat{F}_{\mbox{\footnotesize Lasso}}$ & $\hat{F}_{\mbox{\footnotesize NW}}$ \\
\hline \hline
$F_{\mbox{\footnotesize Linear}}$    & median & {\bf 0.057} & 1.156 & 0.060 & 0.507 \\
                & mean   & 0.095 & 1.124 & {\bf 0.066} & 0.533 \\
                & s.d.   & 0.143 & 0.241 & 0.026 & 0.081 \\
\hline
$F_{\mbox{\footnotesize SI}}$    & median & {\bf 0.050} & 0.502 & 0.795 & 0.308 \\
                & mean   & {\bf 0.051} & 0.539 & 0.776 & 0.326 \\
                & s.d.   & 0.011 & 0.200 & 0.208 & 0.109 \\
\hline
$F_{\mbox{\footnotesize NP}}$    & median & {\bf 0.358} & 0.788 & 1.910 & 0.374 \\
                & mean   & 0.504 & 0.771 & 1.931 & {\bf 0.391} \\
                & s.d.   & 0.320 & 0.168 & 0.468 & 0.101 \\
\hline
\end{tabular}
\vspace*{-6pt}
\begin{tabular}{|p{1.6cm}|p{1.4cm}||p{1.7cm}|p{1.7cm}|p{1.7cm}|p{1.7cm}|p{1.3cm}
|}
\hline $n=100$  &  $p=50$       & $\hat{F}_{\mbox{\footnotesize Fourier}}$ & $\hat{F}_{\mbox{\footnotesize HHI}}$
                             & $\hat{F}_{\mbox{\footnotesize Lasso}}$ & $\hat{F}_{\mbox{\footnotesize NW}}$ \\
\hline \hline
$F_{\mbox{\footnotesize Linear}}$    & median & 0.053 & 0.092 & {\bf 0.050} & 0.519 \\
                & mean   & 0.054  & 0.100 & {\bf 0.050} & 0.508  \\
                & s.d.   & 0.007 & 0.026 & 0.006 & 0.026 \\
\hline
$F_{\mbox{\footnotesize SI}}$    & median & {\bf 0.047} & 0.242 & 0.503 & 0.329 \\
                & mean   & {\bf 0.070} & 0.267 & 0.502 & 0.339 \\
                & s.d.   & 0.099 & 0.111 & 0.106 & 0.073 \\
\hline
$F_{\mbox{\footnotesize NP}}$    & median & {\bf 0.361} & 0.736 & 1.968 & 0.418 \\
                & mean   & 0.557 & 0.765 & 2.045 & {\bf 0.406} \\
                & s.d.   & 0.519 & 0.226 & 0.546 & 0.076 \\
\hline
\hline
\end{tabular}
\caption{\textsf{Numerical results for the simulated data, with $n=50$ and $n=100$,
$p=50$.
The characters in bold indicate the best performance.}}
\label{table1bis}
\end{center}
\end{table}
\medskip

Some comments are in order. First, we note without surprise that:
\begin{enumerate}
 \item The Lasso performs well in the linear setting [{\bf Model 1}].
 \item The single-index methods $\hat F_{\mbox{\footnotesize Fourier}}$ and $\hat F_{\mbox{\footnotesize HHI}}$ are the best ones when the targeted regression function really involves a
single-index model [{\bf Model 2}].
 \item The kernel method gives good results in the purely nonparametric setting [{\bf Model 3}].
\end{enumerate}
Interestingly, $\hat{F}_{\mbox{\footnotesize Fourier}}$ provides slightly better results than the single-index-tailored estimate $\hat{F}_{\mbox{\footnotesize HHI}}$, especially for $p=50$. This observation can be easily explained by the fact that $\hat{F}_{\mbox{\footnotesize HHI}}$ does not integrate any sparsity information regarding the parameter $\theta^{\star}$, whereas $\hat{F}_{\mbox{\footnotesize Fourier}}$ tries to focus on the dimension of the active coordinates, which is equal to 2 in this simulation. As a general finding, we retain that $\hat{F}_{\mbox{\footnotesize Fourier}}$ is the most robust of all the tested procedures. 
\subsection{Real data}
The real-life data sets used in this second series of experiments are from two different sources. The first one, called {\bf AIR-QUALITY} data ($n=111$, $p=3$), has been  first used by Chambers, Cleveland,
Kleiner and Tukey
\cite{chambers} and has been later considered as a benchmark in the study and comparison of single-index models (see, for example, 
Antoniadis, Gr\'egoire and McKeague \cite{Antoniadis} and Wang
\cite{Wang}, among others). This data set originated from an environmental study relating $n=111$ ozone concentration measures at $p=3$ meteorological variables, namely wind speed, temperature and radiation. The data is available as a package in the software $\mathbf R$ \cite{R}, which we employed in all the numerical experiments. The programs are available upon request from the authors.
\medskip

The second category of data arises from the UC Irvine Machine Learning Repository \texttt{http://archive.ics.uci.edu/ml}, where the following packages have been downloaded from: 
\begin{itemize}
\item {\bf AUTO-MPG} (Quinlan \cite{quinlan}, $n=392$, $p=7$).
\item {\bf CONCRETE} (Yeh \cite{concrete}, $n=1030$, $p=8$).
\item {\bf HOUSING} (Harrison and Rubinfeld \cite{housing}, $n=508$, $p=13$). 
\item {\bf SLUMP-1}, {\bf SLUMP-2} and {\bf SLUMP-3}, which correspond to the concrete slump test data introduced by Yeh \cite{concrete2} ($n=51$, $p=7$). Since there are 3 different output variables $Y$ in the original data set, we created a single experiment for each of these variables (1 refers to the output ``slump'', 2 to the output ``flow'' and 3 to the output ``28-day Compressive Strength'').
\item {\bf WINE-RED} and {\bf WINE-WHITE} (Cortez, Cerdeira, Almeida, Ma\-tos and Reis \cite{wine}, $n=1599$, $n=4898$, $p=11$).
\end{itemize}
We refer to the above-mentioned references for a precise description of the meaning of the variables involved in these data sets. For homogeneity reasons, all data were normalized to force the input variables to lie
in $[-1,1]$---in accordance with the setting of our method---and to ensure that all output
variables have standard deviation $0.5$. In two data sets ({\bf AIR-QUALITY} and {\bf AUTO-MPG}) there were some missing values and the corresponding observations were simply removed.
\medskip

For each method and each of the nine data sets, we randomly split the observations in a learning and a test set of equal sizes, computed the estimate on the learning set, evaluated the prediction error on the test set, and repeated this protocol 20 times. The results are summarized in Table \ref{table2}.
\begin{table}[t!]
\begin{center}
\begin{tabular}{|p{3.3cm}|p{1.2cm}||p{1.7cm}|p{1.7cm}|p{1.7cm}|p{1.3cm}|}
\hline 
\hline
Data set  &        & $\hat{F}_{\mbox{\footnotesize Fourier}}$ & $\hat{F}_{\mbox{\footnotesize HHI}}$
                             & $\hat{F}_{\mbox{Lasso}}$ & $\hat{F}_{\mbox{\footnotesize NW}}$ \\
\hline \hline
{\bf AIR QUALITY}    & median & 0.117 & {\bf 0.099} & 0.107 & 0.129 \\
$n=111$                & mean & 0.128 & {\bf 0.096} & 0.113 & 0.130 \\
$p=3$                & s.d.   & 0.044 & 0.029 & 0.029 & 0.035 \\
\hline
{\bf AUTO-MPG}    & median & {\bf 0.044} & 0.049 & 0.070 & 0.068 \\
$n=392$                & mean & 0.051 & {\bf 0.050} & 0.072 & 0.069 \\
$p=7$                & s.d.   & 0.017 & 0.006 & 0.011 & 0.009 \\
\hline
{\bf CONCRETE}    & median & 0.089 & {\bf 0.087} & 0.106 & 0.094 \\
$n=1030$                & mean   & 0.091 & {\bf 0.087} & 0.107 & 0.094 \\
$p=8$                & s.d.   & 0.008 & 0.003 & 0.005 & 0.004 \\
\hline
{\bf HOUSING}    & median & 0.074 & {\bf 0.059} & 0.086 & 0.086 \\
$n=508$                & mean   & 0.076 & {\bf 0.061} & 0.085 & 0.088 \\
$p=11$                & s.d.   & 0.015 & 0.013 & 0.012 & 0.016 \\
\hline
{\bf SLUMP-1}    & median & 0.289 & {\bf 0.171} & 0.201 & 0.208 \\
$n=51$                & mean   & 0.244 & {\bf 0.187} & 0.213 & 0.226 \\
$p=7$                & s.d.   & 0.062 & 0.050 & 0.049 & 0.047 \\
\hline
{\bf SLUMP-2}    & median & 0.219 & 0.196 & {\bf 0.172} & 0.215 \\
$n=51$                & mean   & 0.216 & 0.194 & {\bf 0.171} & 0.213 \\
$p=7$                & s.d.   & 0.053 & 0.025 & 0.019 & 0.022 \\
\hline
{\bf SLUMP-3}    & median & 0.065 & 0.070 & {\bf 0.053} & 0.116 \\
$n=51$                & mean   & 0.073 & 0.079 & {\bf 0.052} & 0.126 \\
$p=7$                & s.d.   & 0.033 & 0.027 & 0.010 & 0.026 \\
\hline
{\bf WINE-RED}    & median & 0.173 & {\bf 0.171} & 0.183 & 0.171 \\
$n=1599$                & mean   & 0.174 & {\bf 0.170} & 0.174 & 0.183 \\
$p=11$                & s.d.   & 0.009 & 0.008 & 0.007 & 0.010 \\
\hline
{\bf WINE-WHITE}    & median & 0.191 & 0.187 & 0.185 & {\bf 0.184} \\
$n=4898$                & mean   & 0.202 & 0.188 & 0.186 & {\bf 0.185} \\
$p=11$                & s.d.   & 0.045 & 0.003 & 0.004 & 0.004 \\
\hline
\hline
\end{tabular}
\end{center}
\caption{\textsf{Numerical results for the real-life data sets. The characters in bold indicate the best performance.}}
\label{table2}
\end{table}
\medskip

\begin{center}
\begin{figure}[!!h]
\centering
\includegraphics*[width=14cm,height=9cm]{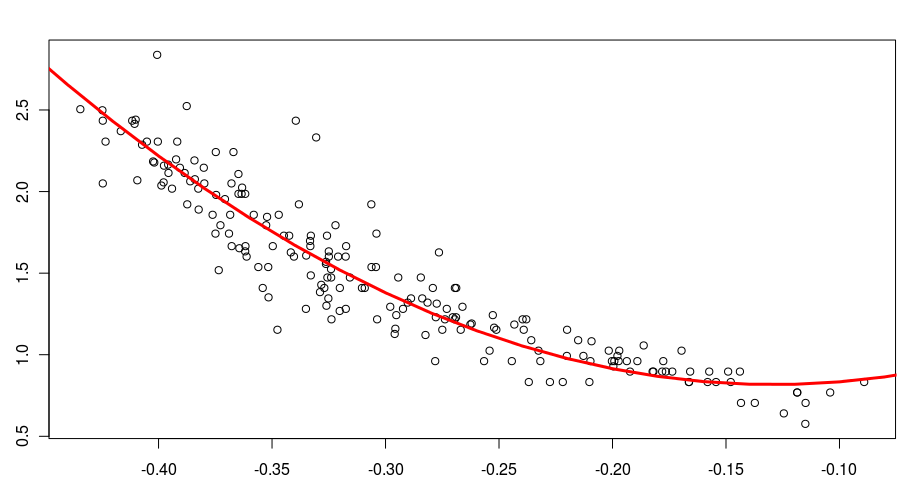}
\caption{\label{fig1} \textsf{{\bf AUTO-MPG} example: Estimated link function by the method $\hat F_{\mbox{\footnotesize Fourier}}$.}}
\end{figure}
\end{center}

We see that all the tested methods provide reasonable results on most data sets. The Lasso is very competitive, especially in the nonasymptotic framework. The estimation procedure $\hat F_{\mbox{\footnotesize Fourier}}$ offers outcomes which are similar to the ones of $\hat{F}_{\mbox{\footnotesize HHI}}$, with a slight advantage for the latter method however. Altogether, $\hat F_{\mbox{\footnotesize Fourier}}$ and $\hat{F}_{\mbox{\footnotesize HHI}}$ provide the best performance in terms of prediction error in 6 out of 9 experiments. Besides, when it is not the best, the method $\hat F_{\mbox{\footnotesize Fourier}}$ is close to the best one, as for example in {\bf SLUMP-3} and {\bf WINE-RED}. As an illustrative example, the plot of the resulting fit of our procedure to the data set {\bf AUTO-MPG} is shown in Figure \ref{fig1}.
Clearly, all data sets under study have a dimension $p$ which is small compared to $n$. To correct this situation, we ran the same series of experiments by adding some additional irrelevant dimensions to the data. Specifically, the observations were embedded into a space of dimension $p\times 4$ by letting the new fake coordinates follow independent uniform $[0,1]$ random variables. The results are shown in Table \ref{table2bis}. In this nonasymptotic framework, the method $\hat{F}_{\mbox{\footnotesize HHI}}$---which is not designed for sparsity---collapses, whereas $\hat F_{\mbox{\footnotesize Fourier}}$ takes a clear advantage over its competitors. In fact, it provides the best results in 3 out of 9 experiments ({\bf AUTO-MPG}, {\bf CONCRETE} and {\bf HOUSING}). Besides, when it is not the best, the method $\hat F_{\mbox{\footnotesize Fourier}}$ is very close to the best one, as for example in {\bf SLUMP-3} and {\bf WINE-RED}. 
\begin{table}[t!]
\begin{center}
\begin{tabular}{|p{3.6cm}|p{1.2cm}||p{1.7cm}|p{1.7cm}|p{1.7cm}|p{1.3cm}|}
\hline 
\hline
Augmented data set  &        & $\hat{F}_{\mbox{\footnotesize Fourier}}$ & $\hat{F}_{\mbox{\footnotesize HHI}}$
                             & $\hat{F}_{\mbox{Lasso}}$ & $\hat{F}_{\mbox{\footnotesize NW}}$ \\
\hline \hline
{\bf AIR QUALITY}    & median & 0.172 & 0.272 & {\bf 0.164} & 0.281 \\
$n=111$                & mean & 0.244 & 0.291 & {\bf 0.163} & 0.291 \\
$p=12$                & s.d.   & 0.163 & 0.116 & 0.038 & 0.046 \\
\hline
{\bf AUTO-MPG}    & median & {\bf 0.043} & 0.062 & 0.085 & 0.202 \\
$n=392$                & mean & {\bf 0.044} & 0.072 & 0.086 & 0.203 \\
$p=28$                & s.d.   & 0.009 & 0.018 & 0.008 & 0.014 \\
\hline
{\bf CONCRETE}    & median & {\bf 0.087} & 0.093 & 0.113 & 0.245 \\
$n=1030$                & mean & {\bf 0.087} & 0.094 & 0.112 & 0.094 \\
$p=32$                & s.d.   & 0.007 & 0.008 & 0.005 & 0.009 \\
\hline
{\bf HOUSING}    & median & {\bf 0.071} & 0.199 & 0.092 & 0.226 \\
$n=508$                & mean   & {\bf 0.075} & 0.181 & 0.095 & 0.227 \\
$p=44$                & s.d.  & 0.023 & 0.084 & 0.013 & 0.018 \\
\hline
{\bf SLUMP-1}    & median & {\bf 0.270} & 0.426 & 0.276 & 0.271 \\
$n=51$                & mean & 0.290 & 0.409 & 0.274 & {\bf 0.262} \\
$p=44$                & s.d. & 0.101 & 0.079 & 0.055 & 0.042 \\
\hline
{\bf SLUMP-2}    & median & 0.276 & 0.332 & {\bf 0.195} & 0.253 \\
$n=51$                & mean & 0.285 & 0.349 & {\bf 0.198} & 0.254 \\
$p=44$                & s.d.  & 0.075 & 0.063 & 0.043 & 0.034 \\
\hline
{\bf SLUMP-3}    & median & 0.079 & 0.371 & {\bf 0.061} & 0.372 \\
$n=51$                & mean  & 0.082 & 0.361 & {\bf 0.058} & 0.279 \\
$p=28$                & s.d.  & 0.025 & 0.079 & 0.013 & 0.031 \\
\hline
{\bf WINE-RED}    & median & 0.178 & 0.222 & {\bf 0.172} & 0.245 \\
$n=1599$                & mean   & 0.176 & 0.226 & {\bf 0.174} & 0.246 \\
$p=44$                & s.d.   & 0.085 & 0.033 & 0.006 & 0.029 \\
\hline
{\bf WINE-WHITE}    & median & 0.199 & 0.239 & {\bf 0.187} & 0.252 \\
$n=4898$                & mean   & 0.204 & 0.256 & {\bf 0.188} & 0.260 \\
$p=11$                & s.d.   & 0.091 & 0.041 & 0.005 & 0.019 \\
\hline
\hline
\end{tabular}
\end{center}
\caption{\textsf{Numerical results for the real-life data sets augmented with noise variables. The characters in bold indicate the best performance.}}
\label{table2bis}
\end{table}
\medskip

Thus, as a general conclusion to this experimental section, we may say that our PAC-Bayesian oriented procedure has an excellent predictive ability, even in nonasymptotic/high-dimensional situations. It is fast, robust, and exhibits performance at the level of the gold standard Lasso.
\section{Proofs}
\label{sectionproofs}
\subsection{Preliminary results}
Throughout this section, we let $\pi$ be the prior probability measure on
$\mathbb R^p \times \mathcal F_n(C+1)$
equipped with its canonical Borel $\sigma$-field.  Recall
that $\mathcal F_n(C+1) \subset \mathcal F$ and that, for each $f \in \mathcal
F_n(C+1)$, we have $\|f\|_{\infty} \leq C+1$.
\medskip

Besides, since $\mathbb E[Y|\bX]=f^{\star}(\theta^{\star T}\bX)$ almost surely, we note once and for all that for all $(\theta,f) \in \mathcal S_{1,+}^p \times \mathcal F_n(C+1)$,
\begin{align*}
R(\theta,f)-R(\theta^{\star},f^{\star})&=\mathbb E \left [Y-f(\theta^T\bX)\right]^2-\mathbb E \left [Y-f^{\star}(\theta^{\star T}\bX)\right]^2\\
&= \mathbb E \left [ f(\theta^T\bX)-f^{\star}(\theta^{\star T}\bX)\right]^2
\end{align*}
(Pythagora's theorem). We start with four technical lemmas. Lemma \ref{lemmemassart} is a
version of Bernstein's inequality, whose proof can be found in Massart
\cite[Chapter 2, inequality (2.21)]{Massart}.  Lemma
\ref{lemmecatoni} is a classical result, whose proof can be found, for example, in Catoni \cite[page 4]{Catoni07}. For a random variable $Z$, the
notation $(Z)_+$ means the positive part of $Z$.
\begin{lem}
\label{lemmemassart}
Let $T_{1}, \hdots, T_{n}$ be independent real-valued random variables. Assume
that there exist two positive constants $v$ and $w$ such that, for all integers
$k\geq 2$,
$$ \sum_{i=1}^{n} \mathbb{E}\left[(T_{i})_{+}^{k}\right] \leq
\frac{k!}{2}vw^{k-2}.$$
Then, for any $\zeta\in\,]0,1/w[$,
$$ \mathbb{E}
\left[\exp\left(\zeta\sum_{i=1}^{n}[T_{i}-\mathbb{E} T_{i}] \right)\right]
        \leq \exp\left(\frac{v\zeta^{2}}{2(1-w\zeta)} \right) .$$
\end{lem}
Given a measurable space $(E,\mathcal E)$ and two probability measures $\mu_1$
and
$\mu_2$ on $(E,\mathcal E)$, we denote by $\mathcal K(\mu_1,\mu_2)$ the
Kullback-Leibler divergence of $\mu_1$ with respect to $\mu_2$, defined by
\begin{equation*}
\mathcal K(\mu_1,\mu_2)=
\left \{ \begin{array}{ll}
\displaystyle \int  \log \left (
\frac{\mbox{d}\mu_1}{\mbox{d}\mu_2}\right)\mbox{d}\mu_1& \mbox{if } \mu_1 \ll
\mu_2,\\
\infty & \mbox{otherwise.} 
\end{array}
\right .
\end{equation*}
(Notation $\mu_1 \ll \mu_2$ means ``$\mu_1$ is absolutely continuous with respect to $\mu_2$''.) 
\begin{lem}
\label{lemmecatoni}
Let $(E,\mathcal E)$ be a measurable space. For any probability measure $\mu$ on
$(E,\mathcal E)$ and any measurable function
$h:E\rightarrow\mathbb{R}$ such that $ \int (\exp \circ h)\emph{d}\mu 
<\infty $,
we have
\begin{equation} \label{lemmacatoni}
\log \int (\exp\circ h)\emph{d}\mu =\sup_{m}\left(\int h\emph{d}m
-\mathcal{K}(m,\mu)\right),
\end{equation}
where the supremum is taken over all probability measures on $(E,\mathcal E)$
and, by convention,
$\infty-\infty=-\infty$. Moreover, as soon as $h$ is bounded from above on the support of $\mu$, the supremum with respect to
$m$ on the right-hand side of (\ref{lemmacatoni}) is reached for the Gibbs
distribution $g$ given by
$$\frac{\emph{d} g}{\emph{d} \mu}(e)=\frac{\exp\left[h(e)\right]}{\displaystyle
\int(\exp \circ h)\emph{d}\mu},  \quad e \in E.$$
\end{lem}
\begin{lem}
\label{thm1}
Assume that Assumption {\bf N} holds. Set  $w= 8(2C+1)\max[L,\newline2C+1]$ and take
$$\lambda \in\,\left]0,
\frac{n}{w+\left[(2C+1)^{2}+4\sigma^{2}\right]}\right[.$$
Then, for all $\delta\in\,]0,1[$ and any data-dependent probability measure
$\hat{\rho}$ absolutely continuous with respect to $\pi$ we have, with
probability at least $1-\delta$,
\begin{align*}
&
              R(\hat\theta,\hat f)-R(\theta^{\star},f^{\star})\\
& \quad  \leq 
       \frac{1}{1-\frac{\lambda
\left[(2C+1)^{2}+4\sigma^{2}\right]}{n-w\lambda}} \Biggr(
          R_{n}(\hat\theta,\hat f) -R_{n}(\theta^{\star},f^{\star})
           +\frac{
\log\left(\frac{\emph{d}\hat{\rho}}{\emph{d}\pi}(\hat\theta,\hat f)\right)+
\log\left(\frac{1}{\delta}\right)}{\lambda}\Biggr)
          ,
\end{align*}
where the pair $(\hat \theta, \hat f)$ is distributed according to $\hat \rho$.
\end{lem}
\noindent{\bf Proof of Lemma \ref{thm1}}. \quad Fix
$\theta\in\mathcal{S}^{p}_{1,+}$ and $f \in \mathcal F_n(C+1)$. The proof starts with an application of Lemma
\ref{lemmemassart} to the random variables
$$ T_{i} =  - \left(Y_{i}-f(\theta^T\bX_{i})\right)^{2}
                    + \left(Y_{i}-f^{\star}(\theta^{\star T}\bX_{i})\right)^{2}
, \quad i=1, \hdots, n.$$
Note that these random variables are independent, identically distributed, and
that
\begin{align*}
&\sum_{i=1}^{n} \mathbb{E} T_{i}^{2}\\
& \quad = \sum_{i=1}^{n} \mathbb{E} \left \{
       \left[2Y_{i} - f(\theta^T\bX_{i})-f^{\star}(\theta^{\star
T}\bX_{i})\right]^{2}
\left[f(\theta^T\bX_{i})-f^{\star}(\theta^{\star T}\bX_{i})\right]^2
            \right\}
\\
&\quad = \sum_{i=1}^{n} \mathbb{E} \left\{
       \left[2\bW_{i} +f^{\star }(\theta^{\star T}\bX_{i}) -
f(\theta^T\bX_{i})\right]^{2}
\left[f(\theta^T\bX_{i})-f^{\star}(\theta^{\star T}\bX_{i})\right]^2
            \right\}
\\
&\quad \leq
\sum_{i=1}^{n} \mathbb{E} \left\{
        \left [4 \bW_{i}^{2} + (2C+1)^{2}\right]
\left[f(\theta^T\bX_{i})-f^{\star}(\theta^{\star T}\bX_{i})\right]^2
            \right\}
\\
& \qquad (\mbox{since } \mathbb E[W_i|\bX_i]=0).
\end{align*}
Thus, by Assumption {\bf N},
\begin{align*}
\sum_{i=1}^{n} \mathbb{E} T_{i}^{2} &\leq
\left[(2C+1)^2+4\sigma^2\right] \sum_{i=1}^{n}
       \mathbb{E}\left[f(\theta^T\bX_{i})-f^{\star}(\theta^{\star
T}\bX_{i})\right]^2
\\
& \leq v, 
\end{align*}
where we set
$$v=  2n [(2C+1)^2+4\sigma^{2}] \left[R(\theta,f)-R(\theta^{\star
},f^{\star})\right].$$
More generally, for all integers $k \geq 3$,
\begin{align*}
&\sum_{i=1}^{n} \mathbb{E}\left[(T_{i})_{+}^{k}\right]\\
& \quad \leq \sum_{i=1}^{n} \mathbb{E} \left\{
       \left|2Y_{i} - f(\theta^T\bX_{i})-f^{\star}(\theta^{\star
T}\bX_{i})\right|^{k}
\left|f(\theta^T\bX_{i})-f^{\star}(\theta^{\star T}\bX_{i})\right|^k
            \right\}
\\
&\quad = \sum_{i=1}^{n} \mathbb{E} \left\{
       \left |2\bW_{i} +f^{\star }(\theta^{\star T}\bX_{i}) -
f(\theta^T\bX_{i})\right|^{k}
\left|f(\theta^T\bX_{i})-f^{\star}(\theta^{\star T}\bX_{i})\right|^k
            \right\}\\
            &\quad \leq 2^{k-1} \sum_{i=1}^{n} \mathbb{E} \left\{
       \left[2^k|\bW_{i}|^{k}
+(2C+1)^{k}\right](2C+1)^{k-2}\left|f(\theta^T\bX_{i})-f^{\star}(\theta^{
\star T}\bX_{i})\right|^{2}
            \right\}.
\end{align*}
In the last inequality, we used the fact that $|a+b|^k\leq 2^{k-1}(|a|^k+|b|^k)$ together with 
\begin{align*}
& \left |f(\theta^T\bX_{i})-f^{\star}(\theta^{
\star T}\bX_{i}) \right|^{k} \\
& \quad =\left |f(\theta^T\bX_{i})-f^{\star}(\theta^{
\star T}\bX_{i}) \right|^{k-2} \times \left |f(\theta^T\bX_{i})-f^{\star}(\theta^{
\star T}\bX_{i}) \right|^{2}\\
& \quad \leq   (2C+1)^{k-2} \left |f(\theta^T\bX_{i})-f^{\star}(\theta^{
\star T}\bX_{i}) \right|^{2}.
\end{align*}
Therefore, by Assumption {\bf N},
\begin{align*}
&\sum_{i=1}^{n} \mathbb{E}\left[(T_{i})_{+}^{k}\right]\\
&\quad \leq \sum_{i=1}^{n} \left[2^{2k-2}{k!}\sigma^{2}L^{k-2}
                            +2^{k-1}(2C+1)^{k}\right] (2C+1)^{k-2}
\left[R(\theta,f)-R(\theta^{\star},f^{\star})\right]
\\
&\quad = v \times \frac{\left[2^{2k-2}{k!}\sigma^{2}L^{k-2}
                            +2^{k-1}(2C+1)^{k}\right] (2C+1)^{k-2}}
       {[(2C+1)^2+4\sigma^{2}]}
\\
&\quad \leq v \times
\frac{8^{k-2} k! \max\left [L^{k-2},(2C+1)^{k-2}\right](2C+1)^{k-2}}{2}\\
&\quad   = \frac{k!}{2}v w^{k-2},
\end{align*}
with $w= 8(2C+1)\max[L,2C+1]$.
\medskip

Thus, for any inverse temperature parameter $\lambda \in\,]0,n/w[$, taking
$\zeta=\lambda/n$, we may write by Lemma \ref{lemmemassart}
\begin{align*}
& \mathbb{E} \Big\{\exp\left[ \lambda
\left(R(\theta,f)-R(\theta^{\star},f^{\star})-R_{n}(\theta,f)+R_{n}(\theta^{
\star},f^{\star})\right)\right]\Big\}\\
&\quad \leq
\exp\left(\frac{v\lambda^{2}}{2n^{2}(1-\frac{w\lambda}{n})}\right).
\end{align*}
Therefore, using the definition of $v$, we obtain
\begin{align*}
&\mathbb{E} \Bigg \{ \exp \Bigg [\left(\lambda
-\frac{\lambda^{2}\left[(2C+1)^{2}+4\sigma^{2}\right]}{n(1-\frac{w\lambda}{n})}
\right)
                      \left(R(\theta,f)-R(\theta^{\star},f^{\star})\right)
                    \\
   &\qquad \qquad  \quad
+\lambda\left(-R_{n}(\theta,f)+R_{n}(\theta^{\star},f^{\star})\right)
         - \log\left(\frac{1}{\delta}\right)\Bigg]\Biggr\}
\leq \delta.
\end{align*}
Next, we use a standard PAC-Bayesian approach (Catoni
\cite{Catoni04,Catoni07}, Audibert \cite{Audibert} and Alquier
\cite{Alquier}). Let us remind the reader that $\pi$ is a prior probability measure on the set $\mathcal S_{1,+}^p\times \mathcal F_n(C+1)$. We have
\begin{align*}
&\int \mathbb{E} \Bigg \{\exp\Biggl[\left(\lambda
-\frac{\lambda^{2}\left[(2C+1)^{2}+4\sigma^{2}\right]}{n(1-\frac{w\lambda}{n})}
\right)
                      \left(R(\theta,f)-R(\theta^{\star},f^{\star})\right)
                    \\
& \qquad \qquad  \qquad      
+\lambda\left(-R_{n}(\theta,f)+R_{n}(\theta^{\star},f^{\star})\right)
         - \log\left(\frac{1}{\delta}\right)\Biggr]\Bigg
\}\mbox{d}\pi(\theta,f)
\leq \delta
\end{align*}
and consequently, using Fubini's theorem,
\begin{align*}
&\mathbb{E} \Bigg\{\int \exp\Biggl[\left(\lambda
-\frac{\lambda^{2}\left[(2C+1)^{2}+4\sigma^{2}\right]}{n(1-\frac{w\lambda}{n})}
\right)
                      \left(R(\theta,f)-R(\theta^{\star},f^{\star})\right)
                    \\
 &\qquad \qquad  \qquad     
+\lambda\left(-R_{n}(\theta,f)+R_{n}(\theta^{\star},f^{\star})\right)-
\log\left(\frac{1}{\delta}\right)\Biggr] \mbox{d}\pi(\theta,f)\Bigg\}
\leq \delta.
\end{align*}
Therefore, for any data-dependent posterior probability measure $\hat{\rho}$
absolutely continuous with respect
to $\pi$, adopting the convention $\infty \times 0=0$,
\begin{align*}
&\mathbb{E}\Bigg\{ \int \exp\Biggl [\left(\lambda
-\frac{\lambda^{2}\left[(2C+1)^{2}+4\sigma^{2}\right]}{n(1-\frac{w\lambda}{n})}
\right)
                      \left(R(\theta,f)-R(\theta^{\star},f^{\star})\right)
                    \\
 &\qquad  
\qquad \qquad
+\lambda\left(-R_{n}(\theta, f)+R_{n}(\theta^{\star},f^{\star})\right)\\
 & \qquad \qquad \qquad
           - \log\left(\frac{\mbox{d}\hat{\rho}}{\mbox{d}\pi}(\theta, f)\right)
         - \log\left(\frac{1}{\delta}\right)\Biggr
]\mbox{d}\hat{\rho}(\theta,f)\Bigg\}\\
&\quad \leq \delta.
\end{align*}
Recalling that $\mathbf P^{\otimes n}$ stands for the distribution of the sample
$\mathcal D_n$, the latter
inequality can be more conveniently written as
\begin{align*}
&\mathbb{E}_{\mathcal D_n \sim\mathbf P^{\otimes n}}
\mathbb{E}_{(\hat{\theta},\hat{f})\sim\hat{\rho}}
   \Biggl \{ \exp\Biggl [\left(\lambda
-\frac{\lambda^{2}\left[(2C+1)^{2}+4\sigma^{2}\right]}{n(1-\frac{w\lambda}{n})}
\right)
                      \left (R(\hat\theta,\hat
f)-R(\theta^{\star},f^{\star})\right)
                    \\
 &\qquad \qquad        +\lambda\left(-R_{n}(\hat \theta,\hat
f)+R_{n}(\theta^{\star},f^{\star})\right)
           - \log\left(\frac{\mbox{d}\hat{\rho}}{\mbox{d}\pi}(\hat\theta,\hat
f)\right)
         - \log\left(\frac{1}{\delta}\right)\Biggr ]\Biggr\}\\
&\quad \leq \delta.
\end{align*}
Thus, using the elementary inequality $\exp(\lambda x) \geq
\mathbf{1}_{\mathbb{R}_{+}}(x)$ we obtain, with probability at most
$\delta$,
\begin{align*}
&
\left(1-\frac{\lambda
\left[(2C+1)^{2}+4\sigma^{2}\right]}{n(1-\frac{w\lambda}{n})}\right)
                      \left(R(\hat\theta,\hat
f)-R(\theta^{\star},f^{\star})\right)
                            \\     
                            &\qquad \qquad \qquad \qquad    \geq
          R_{n}(\hat\theta,\hat f) -R_{n}(\theta^{\star},f^{\star})          
+\frac{\log\left(\frac{\mbox{d}\hat{\rho}}{\mbox{d}\pi}(\hat\theta,\hat
f)\right)+ \log\left(\frac{1}{\delta}\right)}{\lambda},
\end{align*}
where the probability is evaluated with respect to the distribution $\mathbf
P^{\otimes n}$ of the
data $\mathcal D_n$ {\it and} the conditional probability measure $\hat \rho$.
Put differently, letting 
$$\lambda \in\,\left]0,
\frac{n}{w+\left[(2C+1)^{2}+4\sigma^{2}\right]}\right[,$$
we have, with probability at least $1-\delta$,
\begin{align*}
&
              R(\hat\theta,\hat f)-R(\theta^{\star},f^{\star}) \\
& \quad \leq 
       \frac{1}{1-\frac{\lambda
\left[(2C+1)^{2}+4\sigma^{2}\right]}{n-w\lambda}} \Biggr(
          R_{n}(\hat\theta,\hat f) -R_{n}(\theta^{\star},f^{\star})
           +\frac{
\log\left(\frac{\mbox{d}\hat{\rho}}{\mbox{d}\pi}(\hat\theta,\hat f)\right)+
\log\left(\frac{1}{\delta}\right)}{\lambda}\Biggr).
\end{align*}
This concludes the proof of Lemma \ref{thm1}.
\hfill $\blacksquare$
\begin{lem}
\label{thm2BIS}
Under the conditions of Lemma \ref{thm1} we have, with probability at least $1-\delta$,
\begin{align*}
&  \int R_{n}(\theta,f)
\emph{d}\hat{\rho}(\theta,f)-R_{n}(\theta^{\star},f^{\star})  \nonumber \\ 
                      &  \quad  \leq \left(1+\frac{\lambda
\left[(2C+1)^{2}+4\sigma^{2}\right]}{n-w\lambda}\right)
          \left(\int R(\theta,f) \emph{d}\hat{\rho}(\theta,f) -R
(\theta^{\star},f^{\star})\right)\nonumber\\
& \quad \qquad +\frac{
\mathcal{K}\left(\hat{\rho},\pi\right)
         + \log\left(\frac{1}{\delta}\right)}{\lambda}.
\end{align*}
\end{lem}
\noindent{\bf Proof of Lemma \ref{thm2BIS}}. \quad The beginning of the proof is similar to the one of Lemma \ref{thm1}. More precisely, we apply Lemma \ref{lemmemassart} with $T_i=  (Y_{i}-f(\theta^T\bX_{i}))^{2}
                    - (Y_{i}-f^{\star}(\theta^{\star T}\bX_{i}))^{2} $ and obtain, for
any inverse temperature parameter $\lambda \in\,]0,n/w[$,
\begin{align*}
& \mathbb{E} \Big\{\exp\left[ \lambda
\left(R(\theta^\star,f^\star)-R(\theta,f)-R_{n}(\theta^\star,f^\star)+R_{n}(\theta
,f)\right)\right]\Big\}\\
&\quad \leq
\exp\left(\frac{v\lambda^{2}}{2n^{2}(1-\frac{w\lambda}{n})}\right).
\end{align*}
Thus, using the definition of $v$, 
\begin{align*}
&\mathbb{E} \Bigg \{ \exp \Bigg [\left(\lambda
+\frac{\lambda^{2}\left[(2C+1)^{2}+4\sigma^{2}\right]}{n(1-\frac{w\lambda}{n})}
\right)
                      \left(R(\theta^\star,f^\star)-R(\theta,f)\right)
                    \\
   &\qquad \qquad  \quad
+\lambda\left(R_{n}(\theta,f)-R_{n}(\theta^\star,f^\star)\right)
         - \log\left(\frac{1}{\delta}\right)\Bigg]\Biggr\}
\leq \delta.
\end{align*}
Integrating with respect to $\pi$ leads to
\begin{align*}
&\int \mathbb{E} \Bigg \{\exp\Biggl[\left(\lambda
+\frac{\lambda^{2}\left[(2C+1)^{2}+4\sigma^{2}\right]}{n(1-\frac{w\lambda}{n})}
\right)
                      \left(R(\theta^\star,f^\star)-R(\theta,f)\right)
                    \\
& \qquad \qquad  \qquad      
+\lambda\left(R_{n}(\theta,f)-R_{n}(\theta^\star,f^\star)\right)
         - \log\left(\frac{1}{\delta}\right)\Biggr]\Bigg
\}\mbox{d}\pi(\theta,f)
\leq \delta
\end{align*}
whence, by Fubini's theorem,
\begin{align*}
&\mathbb{E} \Bigg\{\int \exp\Biggl[\left(\lambda
+\frac{\lambda^{2}\left[(2C+1)^{2}+4\sigma^{2}\right]}{n(1-\frac{w\lambda}{n})}
\right)
                      \left(R(\theta^\star,f^\star)-R(\theta,f)\right)
                    \\
 &\qquad \qquad  \qquad     
+\lambda\left(R_{n}(\theta,f)-R_{n}(\theta^\star,f^\star)\right)-
\log\left(\frac{1}{\delta}\right)\Biggr] \mbox{d}\pi(\theta,f)\Bigg\}
\leq \delta.
\end{align*}
Thus, for any data-dependent posterior probability measure $\hat{\rho}$
absolutely continuous with respect
to $\pi$,
\begin{align*}
&\mathbb{E}\Bigg\{ \int \exp\Biggl [\left(\lambda
+\frac{\lambda^{2}\left[(2C+1)^{2}+4\sigma^{2}\right]}{n(1-\frac{w\lambda}{n})}
\right)
                      \left(R(\theta^{\star},f^{\star})-R(\theta,f)\right)
                    \\
 &\qquad  
\qquad \qquad
+\lambda\left(R_{n}(\theta, f)-R_{n}(\theta^{\star},f^{\star})\right)\\
 & \qquad \qquad \qquad
           - \log\left(\frac{\mbox{d}\hat{\rho}}{\mbox{d}\pi}(\theta, f)\right)
         - \log\left(\frac{1}{\delta}\right)\Biggr
]\mbox{d}\hat{\rho}(\theta,f)\Bigg\}\\
&\quad \leq \delta.
\end{align*}
Therefore, by Jensen's inequality,
\begin{align*}
&\mathbb{E}\Bigg\{  \exp \int \Biggl [\left(\lambda
+\frac{\lambda^{2}\left[(2C+1)^{2}+4\sigma^{2}\right]}{n(1-\frac{w\lambda}{n})}
\right)
                      \left(R(\theta^{\star},f^{\star})-R(\theta,f)\right)
                    \\
 &\qquad  
\qquad \qquad
+\lambda\left(R_{n}(\theta, f)-R_{n}(\theta^{\star},f^{\star})\right)\\
 & \qquad \qquad \qquad
           - \log\left(\frac{\mbox{d}\hat{\rho}}{\mbox{d}\pi}(\theta, f)\right)
         - \log\left(\frac{1}{\delta}\right)\Biggr
]\mbox{d}\hat{\rho}(\theta,f)\Bigg\}\\
& \quad = \mathbb{E}\Bigg\{  \exp  \Biggl [\left(\lambda
+\frac{\lambda^{2}\left[(2C+1)^{2}+4\sigma^{2}\right]}{n(1-\frac{w\lambda}{n})}
\right)
                      \left(R(\theta^{\star},f^{\star})-
                          \int R(\theta,f)\mbox{d}\hat{\rho}(\theta,f)\right)
                    \\
 &\qquad  
\qquad \qquad
+\lambda\left(\int R_{n}(\theta, f)\mbox{d}\hat{\rho}(\theta,f)
         -R_{n}(\theta^{\star},f^{\star})\right)\\
 & \qquad \qquad \qquad
           - \mathcal{K}\left(\hat{\rho},\pi\right)
         - \log\left(\frac{1}{\delta}\right)\Biggr
]\Bigg\}\\
&\quad \leq \delta.
\end{align*}
Consequently, by the elementary inequality $\exp(\lambda x) \geq
\mathbf{1}_{\mathbb{R}_{+}}(x)$, we obtain, with probability at most
$\delta$,
\begin{align*}
&
\int R_{n}(\theta, f)\mbox{d}\hat{\rho}(\theta,f) -R_{n}(\theta^{\star},f^{\star})    
                            \\     
                            &\quad    \geq
\left(1+\frac{\lambda
\left[(2C+1)^{2}+4\sigma^{2}\right]}{n-w\lambda}\right)
                      \left(\int R(\theta,f)\mbox{d}\hat{\rho}(\theta,f)
-R(\theta^{\star},f^{\star})\right)
 \\               
& \qquad \qquad \qquad \qquad
+\frac{\mathcal{K}\left(\hat{\rho},\pi\right)+ \log\left(\frac{1}{\delta}\right)}{\lambda}.
\end{align*}
Equivalently, with probability at least $1-\delta$,
\begin{align*}
&
\int R_{n}(\theta, f)\mbox{d}\hat{\rho}(\theta,f) -R_{n}(\theta^{\star},f^{\star})    
                            \\     
                            &\quad    \leq
\left(1+\frac{\lambda
\left[(2C+1)^{2}+4\sigma^{2}\right]}{n-w\lambda}\right)
                      \left(\int R(\theta,f)\mbox{d}\hat{\rho}(\theta,f)
-R(\theta^{\star},f^{\star})\right)
 \\               
&\quad \qquad
+\frac{\mathcal{K}\left(\hat{\rho},\pi\right)+ \log\left(\frac{1}{\delta}\right)}{\lambda}.
\end{align*}
\hfill $\blacksquare$
\subsection{Proof of Theorem \ref{thm2}}
The proof starts with an application of Lemma \ref{thm1} with
$\hat{\rho}=\hat{\rho}_{\lambda}$ (the Gibbs distribution) as posterior
distribution.
More precisely, we know that, with probability larger than $1-\delta$,
\begin{align*}
&           
R(\hat{\theta}_{\lambda},\hat{f}_{\lambda})-R(\theta^{\star},f^{\star})\\
& \quad \leq 
       \frac{1}{1-\frac{\lambda
\left[(2C+1)^{2}+4\sigma^{2}\right]}{n-w\lambda}} \Bigg(
          R_{n}(\hat\theta_{\lambda},\hat f_{\lambda})
-R_{n}(\theta^{\star},f^{\star})\\
&\qquad \qquad \qquad \qquad
           +\frac{
\log\left(\frac{\mbox{d}\hat{\rho}_{\lambda}}{\mbox{d}\pi}(\hat{\theta}_{
\lambda},\hat{f}_{\lambda})\right)+
\log\left(\frac{1}{\delta}\right)}{\lambda}\Bigg),
\end{align*}
where the probability is evaluated with respect to the distribution $\mathbf
P^{\otimes n}$ of the
data $\mathcal D_n$ {\it and} the conditional probability measure $\hat
\rho_{\lambda}$. Observe that
\begin{align*}
 \log\left(\frac{\mbox{d}\hat{\rho}_{\lambda}}{\mbox{d}\pi}(\hat\theta_{\lambda
},\hat f_{\lambda})\right)
          &= \log\left(\frac{\exp \left [-\lambda
R_{n}(\hat{\theta}_{\lambda},\hat f_{\lambda})\right]}
                  {\displaystyle \int \exp\left[-\lambda R_{n}(\theta,f)\right]
\mbox{d}\pi(\theta,f) }\right)
\\
&= -\lambda R_{n}(\hat\theta_{\lambda},\hat f_{\lambda}) - \log \int
\exp\left[-\lambda R_{n}(\theta,f)\right] \mbox{d}\pi(\theta,f).
\end{align*}
Consequently, with probability at least $1-\delta$,
\begin{align*}
&              R(\hat\theta_{\lambda},\hat
f_{\lambda})-R(\theta^{\star},f^{\star})
                    \\
& \quad \leq
       \frac{1}{\lambda\left(1-\frac{\lambda
\left[(2C+1)^{2}+4\sigma^{2}\right]}{n-w\lambda}\right)}
           \Bigg(
          -\log \int \exp\left [-\lambda
R_{n}(\theta,f)\right]\mbox{d}\pi(\theta,f)\\
          &\qquad\qquad \qquad\qquad \qquad -\lambda
R_n(\theta^{\star},f^{\star}) + \log\left(\frac{1}{\delta}\right)\Bigg).
          \end{align*}
Next, using Lemma \ref{lemmecatoni} we deduce that, with probability at least
$1-\delta$,
\begin{align*}
  &R(\hat\theta_{\lambda},\hat f_{\lambda})-R(\theta^{\star},f^{\star})\\
                       & \quad  \leq\frac{1}{1-\frac{\lambda
\left[(2C+1)^{2}+4\sigma^{2}\right]}{n-w\lambda}}
                \inf_{\hat{\rho}}\Bigg\{  
          \int R_{n}(\theta,f)\mbox{d}\hat{\rho}(\theta,f)
-R_{n}(\theta^{\star},f^{\star})\\
          &
           \qquad \qquad \qquad \qquad+\frac{ \mathcal{K}(\hat{\rho},\pi)+
\log\left(\frac{1}{\delta}\right)}{\lambda}
            \Bigg\},
\end{align*}
where the infimum is taken over all probability measures on
$\mathcal{S}_{1,+}^{p} \times \mathcal F_n(C+1)$.
In particular, letting $\mathcal{M}(I,M)$ be the set of all probability measures
on $\mathcal{S}_{1,+}^{p}(I)\times
\mathcal F_M(C+1)$, we have, with probability at least $1-\delta$,
\begin{align*}
  &R(\hat\theta_{\lambda},\hat f_{\lambda})-R(\theta^{\star},f^{\star})\\
                       & \quad  \leq\frac{1}{1-\frac{\lambda
\left[(2C+1)^{2}+4\sigma^{2}\right]}{n-w\lambda}}
                \inf_{ \tiny{\begin{array}{c} I\subset\{1,\hdots,p\}
                                 \\ 1\leq M \leq n \end{array}}}
                \inf_{\hat{\rho}\in\mathcal{M}(I,M)}\Bigg\{  
          \int R_{n}(\theta,f)\mbox{d}\hat{\rho}(\theta,f)
-R_{n}(\theta^{\star},f^{\star})\\
          &
           \quad \qquad \qquad+\frac{ \mathcal{K}(\hat{\rho},\pi)+
\log\left(\frac{1}{\delta}\right)}{\lambda}
            \Bigg\}.
\end{align*}
Next, observe that, for $\hat{\rho}\in\mathcal{M}(I,M)$,
\begin{align}
\mathcal{K}(\hat{\rho},\pi) & = \mathcal{K}(\hat{\rho},\mu\otimes\nu)
 = \mathcal{K}(\hat \rho,\mu_{I} \otimes \nu_{M})
                 +
\log\left[\frac{\left(1-\left(\frac{1}{10}\right)^{p}\right)
 \left(1-\left(\frac{1}{10}\right)^{n}\right){p\choose
|I|}}{10^{-|I|-M}}\right] \nonumber\\
& \leq \mathcal{K}(\hat \rho,\mu_{I} \otimes \nu_{M})
 + \log\left[\frac{{p\choose
|I|}}{10^{-|I|-M}}\right].\label{CAC}
\end{align}
Therefore, with probability at least $1-\delta$,
\begin{align}
&R(\hat\theta_{\lambda},\hat f_{\lambda})-R(\theta^{\star},f^{\star})
\nonumber \\
& \quad      \leq \frac{1}{1-\frac{\lambda
\left[(2C+1)^{2}+4\sigma^{2}\right]}{n-w\lambda}}
\inf_{ \tiny{\begin{array}{c} I\subset\{1,\hdots,p\}
                                 \\ 1\leq M \leq n \end{array}}}
      \inf_{\hat \rho\in\mathcal{M}(I,M)}
\Biggl \{
          \int R_{n}(\theta,f)\mbox{d}\hat \rho(\theta,f) 
-R_{n}(\theta^{\star},f^{\star}) \nonumber \\
          & \quad \qquad +\frac{ \mathcal{K}(\hat \rho,\mu_{I} \otimes \nu_{M})
                 +
\log\left[\frac{{p\choose
|I|}}{10^{-|I|-M}}\right] + \log\left(\frac{1}{\delta}\right)}{\lambda}
     \Biggr\}.
\label{one}
\end{align}
By Lemma \ref{thm2BIS} and inequality (\ref{CAC}), for any data-dependent distribution $\hat \rho \in \mathcal M(I,M)$, with probability at least $1-\delta$,
\begin{align}
&  \int R_{n}(\theta,f)
\mbox{d}\hat{\rho}(\theta,f)-R_{n}(\theta^{\star},f^{\star}) \nonumber \\ 
                      &  \quad \leq \left(1+\frac{\lambda
\left[(2C+1)^{2}+4\sigma^{2}\right]}{n-w\lambda}\right)
          \left(\int R(\theta,f) \mbox{d}\hat{\rho}(\theta,f) -R
(\theta^{\star},f^{\star})\right)\nonumber\\
&\quad \qquad+\frac{  \mathcal{K}(\hat \rho,\mu_{I} \otimes \nu_{M}) +
\log\left[\frac{{p\choose
|I|}}{10^{-|I|-M}}\right]+
\log\left(\frac{1}{\delta}\right)}{\lambda}.
\label{two}
\end{align}
Thus, combining inequalities (\ref{one}) and (\ref{two}), we may write, with
probability at least $1-2\delta$,
\begin{align}
&R(\hat\theta_{\lambda},\hat
f_{\lambda})-R(\theta^{\star},f^{\star})\nonumber\\
 &\quad    \leq \frac{1}{1-\frac{\lambda
\left[(2C+1)^{2}+4\sigma^{2}\right]}{n-w\lambda}}
\inf_{ \tiny{\begin{array}{c} I\subset\{1,\hdots,p\}
                                 \\ 1\leq M \leq n \end{array}}}
\inf_{\hat \rho\in\mathcal{M}(I,M)}\Bigg \{ \nonumber\\
&\quad \qquad\left(1+\frac{\lambda
\left[(2C+1)^{2}+4\sigma^{2}\right]}{n-w\lambda}\right)
          \bigg(\int R(\theta,f)\mbox{d}\hat \rho(\theta,f)-R
(\theta^{\star},f^{\star})\bigg)\nonumber\\
          &\quad \qquad +2 \,\frac{
                     \mathcal{K}(\hat \rho,\mu_{I}\otimes\nu_{M})
                 + \log\left [\frac{{p\choose |I|}}{10^{-|I|-M}}\right]
                    + \log\left(\frac{1}{\delta}\right)}{\lambda}
      \Bigg\}.\label{inegalitedegueu}
\end{align}
For any subset $I$ of $\{1, \hdots, p\}$, any positive integer $M\leq n$ and any
$\eta,\gamma\in\,]0,1/{n}]$, let the probability measure
$\rho_{I,M,\eta,\gamma}$ be defined by
$$ \mbox{d}\rho_{I,M,\eta,\gamma}(\theta,f) =
\mbox{d}\rho^{1}_{I,M,\eta}(\theta) \mbox{d}\rho^{2}_{I,M,\gamma}(f),
$$
with
$$
\frac{\mbox{d}\rho^{1}_{I,M,\eta}}{ \mbox{d}\mu_{I}}(\theta)
\propto
\mathbf{1}_{[\|\theta-\theta^{\star}_{I,M}\|_{1}\leq \eta]}
$$
and
$$
\frac{\mbox{d}\rho^{2}_{I,M,\gamma}}{ \mbox{d}\nu_{M}}(f)  \propto
\mathbf{1}_{[\|f-f^{\star}_{I,M}\|_M \leq \gamma]} $$
where, for $f=\sum_{j=1}^{M}\beta_{j}\varphi_{j}\in\mathcal{F}_M(C+1)$, we put
$$ \|f\|_M = \sum_{j=1}^{M}j|\beta_{j}| .$$
With this notation, inequality \eqref{inegalitedegueu} leads to
\begin{align}
&R(\hat\theta_{\lambda},\hat f_{\lambda})-R(\theta^{\star},f^{\star})\nonumber
\\
 &   \leq \frac{1}{1-\frac{\lambda
\left[(2C+1)^{2}+4\sigma^{2}\right]}{n-w\lambda}}
\inf_{ \tiny{\begin{array}{c} I\subset\{1,\hdots,p\}
                                 \\ 1\leq M \leq n
\end{array}}}\inf_{\eta,\gamma>0} \Bigg\{ \nonumber\\        
 &\quad \qquad \left(1+\frac{\lambda
\left[(2C+1)^{2}+4\sigma^{2}\right]}{n-w\lambda}\right)
\Bigg(\int
R(\theta,f)\mbox{d}\rho_{I,M,\eta,\gamma}(\theta,f) -R
(\theta^{\star},f^{\star})\Bigg)\nonumber\\
&\quad \qquad+2\,\frac{\mathcal{K}(\rho_{I,M,\eta,\gamma},\mu_{I}\otimes\nu_{M})
                 + \log\left [\frac{{p\choose |I|}}{10^{-|I|-M}}\right]
                    + \log\left(\frac{1}{\delta}\right)}{\lambda}
      \Biggr\}.\label{inegalitedegueu2}
\end{align}
To finish the proof, we have to control the different terms in
\eqref{inegalitedegueu2}. Note first that
$$ \log {p\choose |I|}
\leq |I|\log \left(\frac{pe}{|I|}\right)
$$
and, consequently,
\begin{equation}
\log\left[\frac{{p\choose |I|}}{10^{-|I|-M}}\right] \leq |I|\log \left( \frac{pe}{|I|}\right) + \left(|I|+ M\right)\log 10 \label{L1}.
\end{equation}
Next,
\begin{align*}
\mathcal{K}(\rho_{I,M,\eta,\gamma},\mu_{I}\otimes\nu_{M})
&=
\mathcal{K}(\rho^{1}_{I,M,\eta}\otimes\rho^{2}_{I,M,\gamma},\mu_{I}\otimes\nu_
{M}) \\
&= \mathcal{K}(\rho^{1}_{I,M,\eta},\mu_{I}) +
\mathcal{K}(\rho^{2}_{I,M,\gamma},\nu_{M})\nonumber .
\end{align*}
By technical Lemma \ref{technique1}, we know that 
$$
\mathcal{K}(\rho^{1}_{I,M,\eta},\mu_{I})
\leq
(|I|-1) \log\left(\max\left[|I|,\frac{4}{\eta}\right]\right).
$$
Similarly, by technical Lemma \ref{technique2},
$$ \mathcal{K}(\rho^{2}_{I,M,\gamma},\nu_{M})  = M \log\left(\frac{C+1}{\gamma}\right).$$
Putting all the pieces together, we are led to
\begin{equation}
\mathcal{K}(\rho_{I,M,\eta,\gamma},\mu_{I}\otimes\nu_{M})
\leq (|I|-1) \log\left(\max\left[|I|,\frac{4}{\eta}\right]\right)+ M \log \left(\frac{C+1}{\gamma}\right). \label{L2}
\end{equation}
Finally, it remains to control the term
$$\int R(\theta,f)\mbox{d}\rho_{I,M,\eta,\gamma}(\theta,f).$$
To this aim, we write
\begin{align*}
&\int R(\theta,f)\mbox{d}\rho_{I,M,\eta,\gamma}(\theta,f)\\
&\quad = \int\mathbb{E} \left[\left(Y-f(\theta^T\bX)\right)^{2}\right]
\mbox{d}\rho_{I,M,\eta,\gamma}(\theta,f)
\\
&\quad =
\int\mathbb{E} \big[\big(Y-f^{\star}_{I,M}(\theta^{\star T}_{I,M}\bX)+
f^{\star}_{I,M}(\theta^{\star T}_{I,M}\bX)-f(\theta^{\star T}_{I,M}\bX)\\
& \quad \qquad +f(\theta^{\star T}_{I,M}\bX)-f(\theta^T\bX)\big)^{2}\big]
\mbox{d}\rho_{I,M,\eta,\gamma}(\theta,f)\\
&\quad = R(\theta^{\star}_{I,M},f^{\star}_{I,M})\\
& \qquad  + \int\mathbb{E} \Bigl[
\left(f^{\star}_{I,M}(\theta^{\star T}_{I,M}\bX)-f(\theta^{\star
T}_{I,M}\bX)\right)^{2}\\
& \qquad + \left(f(\theta^{\star T}_{I,M}\bX)-f(\theta^T\bX)\right)^{2}
\\
&\qquad  + 2 \left(Y-f^{\star}_{I,M}(\theta^{\star T}_{I,M}\bX)\right)
\left(f^{\star}_{I,M}(\theta^{\star T}_{I,M}\bX)-f(\theta^{\star
T}_{I,M}\bX)\right)
\\
&\qquad +  2\left(Y-f^{\star}_{I,M}(\theta^{\star
T}_{I,M}\bX)\right)\left(f(\theta^{\star T}_{I,M}\bX)-f(\theta^T\bX)\right)
\\
&\qquad + 2 \left(f^{\star}_{I,M}(\theta^{\star T}_{I,M}\bX)-f(\theta^{\star
T}_{I,M}\bX)\right)\left(f(\theta^{\star T}_{I,M}\bX)-f(\theta^T\bX)\right)
\Bigr] \mbox{d}\rho_{I,M,\eta,\gamma}(\theta,f)
\\
&\quad := R(\theta^{\star}_{I,M},f^{\star}_{I,M}) + \mathbf A+\mathbf B+\mathbf
C+\mathbf D+\mathbf E.
\end{align*}
\paragraph{Computation of C}
By Fubini's theorem,
\begin{align*}
{\bf C}  &= \mathbb{E}\left[\int 2 \left(Y-f^{\star}_{I,M}(\theta^{\star
T}_{I,M}\bX)\right)
\left(f^{\star}_{I,M}(\theta^{\star T}_{I,M}\bX)-f(\theta^{\star
T}_{I,M}\bX)\right) \mbox{d}\rho_{I,M,\eta,\gamma}(\theta,f) \right]
\\
&= \mathbb{E}\Biggl\{\int \Biggl[2 \left(Y-f^{\star}_{I,M}(\theta^{\star
T}_{I,M}\bX)\right)
\\
&\quad \qquad \qquad \times \int
\left(f^{\star}_{I,M}(\theta^{\star T}_{I,M}\bX)-f(\theta^{\star
T}_{I,M}\bX)\right)
\mbox{d}\rho^{2}_{I,M,\gamma}(f)
\Biggr]
\mbox{d}\rho^{1}_{I,M,\eta}(\theta)
\Biggr\}.
\end{align*}
By the triangle inequality, for $f=\sum_{j=1}^M \beta_j \varphi_j$ and $f^{\star}_{I,M}=\sum_{j=1}^M (\beta^{\star}_{I,M})_j \varphi_j$, it holds
$$\sum_{j=1}^M j |\beta_j|\leq \sum_{j=1}^M j \left | \beta_j-(\beta^{\star}_{I,M})_j \right|+ \sum_{j=1}^M j \left |(\beta^{\star}_{I,M})_j\right|.$$
Since $f^{\star}_{I,M} \in \mathcal F_M(C)$, we have $\sum_{j=1}^M j | (\beta^{\star}_{I,M})_j| \leq C$,
so that $\sum_{j=1}^M j |\beta_j|\leq C+1$ as soon as $\| f-f^{\star}_{I,M}\|_M \leq 1$. This shows that the set
$$\left\{f = \sum_{j=1}^M \beta_j \varphi_j : \| f-f^{\star}_{I,M}\|_M \leq \gamma\right\}$$
is contained in the support of $\nu_M$. In particular, this implies that $\rho^{2}_{I,M,\gamma}$ is centered at $f^{\star}_{I,M}$ and, consequently, 
$$
\int \left(f^{\star}_{I,M}(\theta^{\star
T}_{I,M}\bX)-f(\theta^{\star T}_{I,M}\bX)\right)
\mbox{d}\rho^{2}_{I,M,\gamma}(f) = 0.$$
This proves
that $\mathbf C=0$. 
\paragraph{Control of A}
Clearly,
$$
{\bf A} \leq  \int \sup_{y\in\mathbb{R}} \left((f^{\star}_{I,M}(y)-f(y)\right)^{2}
\mbox{d}\rho^{2}_{I,M,\gamma}(f)
\leq \gamma^{2}.
$$
\paragraph{Control of B}
We have
\begin{align*}
{\bf B} &= \int\mathbb{E} \left[\left(f(\theta^{\star
T}_{I,M}\bX)-f(\theta^T\bX)\right)^{2} \right]
\mbox{d}\rho_{I,M,\eta,\gamma}(\theta,f)
\\
&\leq \int \mathbb{E} \left[\left(\ell(C+1) (\theta^{\star
T}_{I,M}-\theta^T)\bX\right)^2\right] \mbox{d}\rho^{1}_{I,M,\eta}(\theta)
\\
& \quad (\mbox{by the mean value theorem})\\
&\leq \ell^2(C+1)^{2}\mathbb{E}\left [\|\bX\|^2_{\infty}\right] \int
\|\theta^{\star}_{I,M}-\theta\|_{1}^{2} \mbox{d}\rho^{1}_{I,M,\eta}(\theta)
\\
   & \leq  \ell^2(C+1)^{2}  \eta^{2}\\
   &\quad (\mbox{by Assumption } \mathbf D).
\end{align*}
\paragraph{Control of E}
Write
\begin{align*}
|{\bf E}| &\leq 2 \int\mathbb{E}
\Bigl[
\left|f^{\star}_{I,M}(\theta^{\star T}_{I,M}\bX)-f(\theta^{\star
T}_{I,M}\bX)\right|
\\
& \qquad \qquad \times \left|f(\theta^{\star
T}_{I,M}\bX)-f(\theta^T\bX)\right|
\Bigr]
\mbox{d}\rho_{I,M,\eta,\gamma}(\theta,f)
\\
&\leq
2 \int\mathbb{E}
\Bigl[
\left|f^{\star}_{I,M}(\theta^{\star T}_{I,M}\bX)-f(\theta^{\star
T}_{I,M}\bX)\right|
\\
&\qquad \qquad \times \ell(C+1)\left|(\theta^{\star
T}_{I,M}-\theta^T)\bX\right|
\Bigr]
\mbox{d}\rho_{I,M,\eta,\gamma}(\theta,f)
\\
&\leq 2 \left(\int\mathbb{E} \left[\left(f^{\star}_{I,M}(\theta^{\star
T}_{I,M}\bX)-f(\theta^{\star T}_{I,M}\bX)\right)^{2}\right]
\mbox{d}\rho_{I,M,\eta,\gamma}(\theta,f)\right)^{\frac{1}{2}}
\\
&\quad \qquad \left(\int\mathbb{E}\left[\left(\ell(C+1) (\theta^{\star
T}_{I,M}-\theta^T)\bX\right)^{2}\right]
\mbox{d}\rho_{I,M,\eta,\gamma}(\theta,f)\right)^{\frac{1}{2}}
\\
&\quad (\mbox{by the Cauchy-Schwarz inequality})\\
&\leq 2 \left(\gamma^{2}\right)^{\frac{1}{2}}
\left(\ell^2(C+1)^{2} \eta^{2}\right)^{\frac{1}{2}}\\
&= 2 \ell(C+1) \gamma\eta.
\end{align*}
\paragraph{Control of D}
Finally,\\
\begin{align*}
{\bf D} &= 2 \int\mathbb{E}
\left[
\left(Y-f^{\star}_{I,M}(\theta^{\star T}_{I,M}\bX)\right)\left(f(\theta^{\star
T}_{I,M}\bX)-f(\theta^T\bX)\right)
\right]
\mbox{d}\rho_{I,M,\eta,\gamma}(\theta,f)
\\
&= 2 \int\mathbb{E}
\left[
\left(Y-f^{\star}_{I,M}(\theta^{\star
T}_{I,M}\bX)\right)\left(f^{\star}_{I,M}(\theta^{\star
T}_{I,M}\bX)-f^{\star}_{I,M}(\theta^T\bX)\right)
\right]
\mbox{d}\rho^{1}_{I,M,\eta}(\theta)
\\
& \quad (\mbox{since }\int f\mbox{d}\rho^{2}_{I,M,\gamma}(f)=f^{\star}_{I,M})\\
&= 2 \mathbb{E}
\left[
\left(Y-f^{\star}_{I,M}(\theta^{\star
T}_{I,M}\bX)\right)\int\left(f^{\star}_{I,M}(\theta^{\star
T}_{I,M}\bX)-f^{\star}_{I,M}(\theta^T\bX)\right)
\mbox{d}\rho^{1}_{I,M,\eta}(\theta)
\right]
\\
&\leq
2 \sqrt{\mathbb{E}\left[
\left(Y-f^{\star}_{I,M}(\theta^{\star T}_{I,M}\bX)\right)^{2}\right]}
\\
&\quad \times \sqrt{\mathbb{E}\left[
\int\left(f^{\star}_{I,M}(\theta^{\star
T}_{I,M}\bX)-f^{\star}_{I,M}(\theta^T\bX)\right)
\mbox{d}\rho^{1}_{I,M,\eta}(\theta)
\right]^{2}}
\\
& \quad (\mbox{by the Cauchy-Schwarz inequality})\\
&= 2 \sqrt{R(\theta^{\star}_{I,M},f^{\star}_{I,M})} \sqrt{\mathbb{E}\left[
\int\left(f^{\star}_{I,M}(\theta^{\star
T}_{I,M}\bX)-f^{\star}_{I,M}(\theta^T\bX)\right)
\mbox{d}\rho^{1}_{I,M,\eta}(\theta)
\right]^{2}}.
\end{align*}
The inequality
\begin{align*}
\left | f^{\star}_{I,M}(\theta^{\star T}_{I,M}\bX)-f^{\star}_{I,M}(\theta^T\bX) \right|
&\leq  \ell (C+1) \left |(\theta^{\star T}_{I,M}-\theta^T)\bX\right|\\
& \leq \ell (C+1) \|\theta^{\star}_{I,M}-\theta\|_1
\end{align*}
leads to
\begin{align*}
&\left [\int \left(f^{\star}_{I,M}(\theta^{\star
T}_{I,M}\bX)-f^{\star}_{I,M}(\theta^T\bX)\right)\mbox{d}\rho^{1}_{I,M,\eta}(
\theta)\right]^2\\
& \quad \leq \ell ^2 (C+1)^2 \left [ \int  \|\theta^{\star}_{I,M}-\theta\|_1 \mbox{d}\rho^{1}_{I,M,\eta}(\theta)\right]^2.
\end{align*}
Consequently,
$$
\left[\int \left(f^{\star}_{I,M}(\theta^{\star
T}_{I,M}\bX)-f^{\star}_{I,M}(\theta^T\bX)\right)\mbox{d}\rho^{1}_{I,M,\eta}(
\theta)\right]^{2}
\leq \ell^2(C+1)^2\eta^2,
$$
and therefore
\begin{align*}
\mathbf D &\leq 2 \ell (C+1)\eta \sqrt{R(0,0)/2}\\
 & \leq \sqrt 2 \ell (C+1) \eta\sqrt{C^{2}+\sigma^{2}}.
\end{align*}
Thus, taking $\eta=\gamma=1/n$ and putting all the pieces together, we obtain
$$ \mathbf A+\mathbf B+\mathbf C+\mathbf D+\mathbf E \leq \frac{\Xi_{1}}{n},$$
where
$ \Xi_{1}$ is a positive constant, function of $C$, $\sigma$ and $\ell$. Combining this inequality with \eqref{inegalitedegueu2}-\eqref{L2} yields,
with probability larger than $1-2\delta$,
\begin{align*}
&R(\hat\theta_{\lambda},\hat f_{\lambda})-R(\theta^{\star},f^{\star})\\
&  \quad \leq \frac{1}{1-\frac{\lambda
\left[(2C+1)^{2}+4\sigma^{2}\right]}{n-w\lambda}}
\inf_{ \tiny{\begin{array}{c} I\subset\{1,\hdots,p\}
                                 \\ 1\leq M \leq n \end{array}}}
\Biggl\{ \left(1+\frac{\lambda
\left[(2C+1)^{2}+4\sigma^{2}\right]}{n-w\lambda}\right)
          \Bigg( R(\theta^{\star}_{I,M},f^{\star}_{I,M})\\
          &\quad \qquad -R (\theta^{\star},f^{\star})+ \frac{
\Xi_{1}}{n}\Bigg) +2\,\frac{ M\log(10(C+1)n)
                 + |I|\log(40epn)
                    + \log\left(\frac{1}{\delta}\right)}{\lambda}
      \Biggr\}.
\end{align*}
Choosing finally 
$$\lambda = \frac{n}{w+2\left[(2C+1)^{2}+4\sigma^{2}\right]},$$
we obtain that there exists a positive constant $\Xi_2$, function of $L$,
$C$, $\sigma$ and $\ell$ such that,
with probability at least $1-2\delta$,
\begin{align*}
R(\hat\theta_{\lambda},\hat f_{\lambda})-R(\theta^{\star},f^{\star})      
&\leq \Xi_2
\inf_{ \tiny{\begin{array}{c} I\subset\{1,\hdots,p\}
                                 \\ 1\leq M \leq n \end{array}}}
\Biggl\{
           R(\theta^{\star}_{I,M},f^{\star}_{I,M}) -R
(\theta^{\star},f^{\star})
\\
         &\quad \qquad   + \frac{ M\log(10Cn)
                 + |I|\log(40epn)
                    + \log\left(\frac{1}{\delta}\right)}{n}
      \Biggr\}.
\end{align*}
This concludes the proof of Theorem \ref{thm2}.
\subsection{Proof of Corollary \ref{thm3}}
We already know, by Theorem \ref{thm2}, that with probability
at least $1-\delta$,
\begin{align*}
R(\hat{\theta}_{\lambda},\hat{f}_{\lambda})- R(\theta^{\star},f^{\star}) &\leq
\Xi\inf_{
\tiny{\begin{array}{c} I\subset\{1,\hdots,p\}
                                 \\ 1\leq M \leq n \end{array}}} \Bigg\{
           R(\theta^{\star}_{I,M},f^{\star}_{I,M}) -
R(\theta^{\star},f^{\star}) \\
&\qquad \qquad  + \frac{ M\log(Cn) + |I|\log({pn})+
\log\left(\frac{2}{\delta}\right)}{n}\Bigg\}.
\end{align*}
By definition, for all
$(\theta,f)\in\mathcal{S}_{1,+}^{p}(I)\times\mathcal{F}_{M}(C)$,
$$ R(\theta^{\star}_{I,M},f^{\star}_{I,M}) \leq R(\theta,f) .$$
In particular, if $I^\star$ is such that
$\theta^{\star}\in\mathcal{S}_{1,+}^{p}(I^{\star})$,
then
\begin{align}
\label{proofcor1}
R(\hat{\theta}_{\lambda},\hat{f}_{\lambda})- R(\theta^{\star},f^{\star}) &\leq
\Xi\inf_{\tiny{\begin{array}{c}  1\leq M \leq n \\
f\in\mathcal{F}_{M}(C)
\end{array}}} \Bigg\{
           R(\theta^{\star},f) -
R(\theta^{\star},f^{\star}) \nonumber\\
&\qquad \qquad  + \frac{ M\log(Cn) +
|I^\star|\log({pn})+
\log\left(\frac{2}{\delta}\right)}{n}\Bigg\}.
\end{align}
Observe that, for any $f\in\mathcal{F}_{M}(C)$,
\begin{align*}
R(\theta^{\star},f) -
R(\theta^{\star},f^{\star}) &= \int_{\mathbb{R}^{p}}
\left[f\left(\theta^{\star
T}\bx\right)-f^{\star}\left(\theta^{\star T}\bx
\right)\right]^{2} \mbox{d}\mathbf{P}(\bx,y)
\\
&\leq
B^{2} \int_{-1}^{1} \left[f\left(t\right)-f^{\star}\left(t\right)\right]^{2}
\mbox{d}t.
\end{align*}
Since $f^\star \in L_2\left ([-1,1]\right)$, we may write
$$ f^{\star} = \sum_{j=1}^{\infty}\beta_j^{\star}\varphi_{j} $$
and apply \eqref{proofcor1} with
$$ f  = \sum_{j=1}^{M}\beta^{\star}_j\varphi_{j}. $$
In order to do so, we just need to check that
$f\in\mathcal{F}_{M}(C)$, that is
$$ \sum_{j=1}^{M}j |\beta^{\star}_{j}| \leq C.$$
But, by the Cauchy-Schwarz inequality,
\begin{align*}
\sum_{j=1}^{M}j |\beta^{\star}_{j}| &=\sum_{j=1}^{M}j^{k}
|\beta^{\star}_{j}| j^{1-k}\\
& \leq \sqrt{ \sum_{j=1}^{M}j^{2k}
(\beta^{\star}_{j})^{2} } \sqrt{\sum_{j=1}^M j^{2-2k}}.
\end{align*}
Thus,
\begin{align*}
\sum_{j=1}^{M}j |\beta^{\star}_{j}|& \leq \frac{\pi}{\sqrt 6}\sqrt{ \sum_{j=1}^{M}j^{2k}
(\beta^{\star}_{j})^{2} }\\
&  \quad (\mbox{since, by assumption, } k\geq 2)\\
& \leq C\\
& \quad (\mbox{since } f^{\star} \in \mathcal W(k,{6}C^2/\pi^2)).
\end{align*}
Next, with this choice of $f$,
$$ \int_{-1}^{1} \left[f\left(t\right)-f^{\star}\left(t\right)\right]^{2}
\mbox{d}t
\leq
            \Lambda M^{-2k} $$
for some positive constant $\Lambda$ depending only on $k$ and $C$
(see for instance Tsybakov \cite{Tsybakov}). Therefore, inequality
\eqref{proofcor1} leads to
\begin{align}
R(\hat{\theta}_{\lambda},\hat{f}_{\lambda})- R(\theta^{\star},f^{\star}) &\leq
\Xi\inf_{1\leq M \leq n} \Bigg\{
           \Lambda M^{-2k} \nonumber \\
&\qquad \qquad  + \frac{ M\log(Cn) +
|I^\star|\log({pn})+
\log\left(\frac{2}{\delta}\right)}{n}\Bigg\}.\label{proofcor2}
\end{align}
Letting $\lceil . \rceil$ be the ceiling function and choosing $M=\lceil(n/\log(Cn))^{\frac{1}{2\beta+1}}\rceil$ in
(\ref{proofcor2}) concludes
the proof. 
 \subsection{Some technical lemmas}
\begin{lem}
\label{technique1}
For any subset $I$ of $\{1, \hdots, p\}$, any positive integer $M\leq n$ and any
$\eta\in\,]0,1/{n}]$, let the probability measure
$\rho_{I,M,\eta}^1$ be defined by
$$
\frac{\emph{d}\rho^{1}_{I,M,\eta}}{ \emph{d}\mu_{I}}(\theta)
\propto
\mathbf{1}_{[\|\theta-\theta^{\star}_{I,M}\|_{1}\leq \eta]}.
$$
Then
$$
\mathcal{K}(\rho^{1}_{I,M,\eta},\mu_{I})
\leq
(|I|-1) \log\left(\max\left[|I|,\frac{4}{\eta}\right]\right).
$$
\end{lem}
\noindent{\bf Proof}.\quad For simplicity, we assume that $I=\{1,\hdots,|I|\}$. Up to a permutation of the coordinates, the proof remains valid for any subset $I$ of $\{1, \hdots, p\}$. Still for simplicity, we let $\tilde{\theta}$ denote
$\theta^{\star}_{I,M}$. By a symmetry argument, it can be assumed that $\tilde{\theta}$ has
nonnegative coordinates---this just means that $\tilde{\theta}$ is arbitrarily fixed in one of the
$2^{|I|-1}$ faces of $\mathcal{S}_{1,+}^{p}(I)$. We denote by $\mathcal{FA}$ this face and note that
$$ \mathcal{FA} = \left\{ \theta \in(\mathbb{R}_{+})^{|I|} \times \{0\}^{p-|I|}:
                    \sum_{j=1}^{|I|} \theta_{j} = 1   \right\} .$$
Finally, without loss of generality, we suppose that the largest coordinate in $\tilde{\theta}$
is $\tilde{\theta}_{1}$, and let $\nu$ be the uniform probability measure on $\mathcal{FA}$, defined by
$$\frac{\mbox{d}\nu}{\mbox{d}\mu_{I}}(\theta) = 2^{|I|-1} \mathbf{1}_{[\theta\in\mathcal{FA}]}. $$
\medskip

Set $u=\min(1/|I|,\eta/2)$, and let
\begin{equation*}
\begin{array}{ccc}
 T_2 & = (\tilde{\theta}_{1} - u,  \tilde{\theta}_{2} + u,\tilde{\theta}_{3},
 \hdots, \tilde{\theta}_{|I|} , 0, \dots ,0 ), \\
 T_3 & = (\tilde{\theta}_{1} - u,  \tilde{\theta}_{2} ,\tilde{\theta}_{3} + u,
 \hdots, \tilde{\theta}_{|I|} , 0, \dots ,0 ), \\
 \vdots & \qquad \qquad \vdots \\
 T_{|I|}& =(\tilde{\theta}_{1} - u,  \tilde{\theta}_{2} , \tilde{\theta}_{3} ,
 \hdots, \tilde{\theta}_{|I|} + u, 0, \dots ,0 ).
\end{array}
\end{equation*}
Note that $u\leq 1/|I| \leq \tilde{\theta}_{1}$. Therefore, for each $j$, all the coordinates
of $T_{j}$ are nonnegative. Obviously $\|T_{j}\|_{1}=1$, so that, for all $j$,
$T_j\in \mathcal{FA}$. Denoting by $K$ the convex hull of the set $\{\tilde{\theta},T_2,\hdots,T_{|I|}\}$, we also have $K\subset \mathcal{FA}$. Next, observe that $\|T_j-\tilde{\theta}\|_{1} = 2 u \leq
\eta$, which implies $K \subset \{\theta \in \mathbb R^p: \|\theta-\tilde{\theta}\|_{1}\leq \eta\}$.
\medskip

Clearly,
\begin{align*}
\mathcal{K}(\rho^{1}_{I,M,\eta},\mu_{I})
& = \log \left (\frac{1}{\displaystyle \int \mathbf{1}_{[\|\theta-\theta^{\star}_{I,M}\|_{1}\leq \eta]}
             \mbox{d}\mu_{I}(\theta)}\right) \\
& \leq \log \left (\frac{1}{\displaystyle \int \mathbf{1}_{[\theta\in\mathcal{FA}]}
            \mathbf{1}_{[\|\theta-\theta^{\star}_{I,M}\|_{1}\leq \eta]}
             \mbox{d}\mu_{I}(\theta)}\right).
\end{align*}
Thus, \begin{align*}
\mathcal{K}(\rho^{1}_{I,M,\eta},\mu_{I})
& \leq \log \left (\frac{2^{|I|-1}}{\displaystyle \int
            \mathbf{1}_{[\|\theta-\theta^{\star}_{I,M}\|_{1}\leq \eta]}
             \mbox{d}\nu(\theta)}\right) \\
& \leq \log \left (\frac{2^{|I|-1}}{\displaystyle \int
            \mathbf{1}_{[\theta\in K]}
             \mbox{d}\nu(\theta)}\right).
\end{align*}
Observe that $K$ is homothetic to $\mathcal{FA}$, by a factor of $u$. This means that
$$ \int \mathbf{1}_{[\theta\in K]}
             \mbox{d}\nu(\theta) = u^{|I|-1} . $$
Consequently, we obtain
$$
\mathcal{K}(\rho^{1}_{I,M,\eta},\mu_{I}) \leq \log\left( \left(\frac{2}{u}\right)^{|I|-1} \right)
      \leq (|I|-1) \log\left(\max\left[|I|,\frac{4}{\eta}\right]\right).
$$
\hfill $\blacksquare$
\begin{lem}
\label{technique2}
For any subset $I$ of $\{1, \hdots, p\}$, any positive integer $M\leq n$ and any
$\gamma\in\,]0,1/{n}]$, let the probability measure $\rho^{2}_{I,M,\gamma}$ be defined by
$$
\frac{\emph{d}\rho^{2}_{I,M,\gamma}}{ \emph{d}\nu_{M}}(f)  \propto
\mathbf{1}_{[\|f-f^{\star}_{I,M}\|_M \leq \gamma]} $$
where, for $f=\sum_{j=1}^{M}\beta_{j}\varphi_{j}\in\mathcal{F}_M(C+1)$, we put
$$ \|f\|_M = \sum_{j=1}^{M}j|\beta_{j}| .$$
Then
$$
 \mathcal{K}(\rho^{2}_{I,M,\gamma},\nu_{M}) = M \log\left(\frac{C+1}{\gamma}\right).$$
\end{lem}
\noindent{\bf Proof}.\quad Observe that
$$
 \mathcal{K}(\rho^{2}_{I,M,\gamma},\nu_{M})=\int \log \left (\frac{\mbox{d} \rho^2_{I,M,\gamma}}{\mbox{d} \nu_M}(f) \right)\mbox{d}\rho^2_{I,M,\gamma}(f).$$
Now, 
$$\frac{\mbox{d} \rho^2_{I,M,\gamma}}{\mbox{d} \nu_M}(f)= \frac{\mathbf 1_{[\|f-f^{\star}_{I,M}\|_M \leq \gamma]}(f)}{\zeta},$$
where $\zeta=\int \mathbf 1_{[\|f-f^{\star}_{I,M}\|_M \leq \gamma]}(f)\mbox{d}\nu_M(f).$
It easily follows, using the fact that the support of $\rho^2_{I,M,\gamma}$ is included in the set $\{f \in \mathcal F_M(C+1): \|f-f^{\star}_{I,M}\|\leq \gamma \}$, that
$$ \mathcal{K}(\rho^{2}_{I,M,\gamma},\nu_{M})=\log (1/\zeta).$$
Note that
\begin{align*}
\zeta & = \int \mathbf 1_{[\|f-f^{\star}_{I,M}\|_M \leq \gamma]}(f)\mbox{d}\nu_M(f)\\
& =\frac{\displaystyle \int \mathbf 1_{[\sum_{j=1}^M j |\beta_j-(\beta^{\star}_{I,M})_j|\leq \gamma]}(\beta)\mathbf 1_{[\sum_{j=1}^M j |\beta_j|\leq C+1]}(\beta)\mbox{d}\beta}{\displaystyle \int\mathbf 1_{[\sum_{j=1}^M j |\beta_j|\leq C+1]}(\beta)\mbox{d}\beta},
\end{align*}
where the second equality is true since $\nu_M$ is (the image of) the uniform probability measure on $\{\beta \in \mathbb R^M : \sum_{j=1}^M j |\beta_j| \leq C+1\}$.  This implies
$$
 \mathcal{K}(\rho^{2}_{I,M,\gamma},\nu_{M})=\log \left ( \frac{\displaystyle \int \mathbf 1_{[ \sum_{j=1}^M j |\beta_j|\leq C+1]} (\beta)\mbox{d}\beta}{\displaystyle \int \mathbf 1_{[\sum_{j=1}^M j | \beta_j-(\beta^{\star}_{I,M})_j|\leq \gamma]}(\beta)\mathbf 1_{[\sum_{j=1}^M j |\beta_j|\leq C+1]}(\beta)\mbox{d}\beta}\right).
 $$
 By the triangle inequality,
 $$\sum_{j=1}^M j |\beta_j|\leq \sum_{j=1}^M j \left | \beta_j-(\beta^{\star}_{I,M})_j \right|+ \sum_{j=1}^M j \left |(\beta^{\star}_{I,M})_j\right|.$$
Since $f^{\star}_{I,M} \in \mathcal F_M(C)$, we have $\sum_{j=1}^M j | (\beta^{\star}_{I,M})_j| \leq C$,
so that
$$\mathbf 1_{[\sum_{j=1}^M j |\beta_j|\leq C+1]} \geq  \mathbf 1_{[\sum_{j=1}^M j | \beta_j-(\beta^{\star}_{I,M})_j|\leq \gamma]}$$
as soon as $\gamma\leq 1$. We conclude that
$$
 \mathcal{K}(\rho^{2}_{I,M,\gamma},\nu_{M}) = \log \left(
    \frac{\displaystyle \int \mathbf{1}_{[\sum_{j=1}^{M} j|\beta_{j}|\leq C+1]} \mbox{d}\beta }
  {\displaystyle \int \mathbf{1}_{[\sum_{j=1}^{M} j|\beta_{j}-(\beta_{I,M}^{\star})_{j}|\leq\gamma]} \mbox{d}\beta }\right)
 = M \log\left(\frac{C+1}{\gamma}\right).$$
 \hfill $\blacksquare$
\section{Annex: Description of the MCMC algorithm}
\label{sectionmcmc}
This annex is intended to make thoroughly clear the specification of the proposal conditional densities $k_{1}$ and
$k_{2}$ introduced in Section \ref{sectionexpe}.
\subsection{Notation}
To provide explicit formulas for the conditional densities $k_{1}((\tau,h)|(\theta,f))$ and $k_{2}((\tau,h)|(\theta,f)),$
we first set
$$ f = \sum_{j=1}^{m_{f}} \beta_{f,j} \varphi_{j} \quad \mbox{and} \quad h = \sum_{j=1}^{m_{h}} \beta_{h,j} \varphi_{j},
$$
where it is recalled that $\{\varphi_{j}\}_{j=1}^{\infty}$ denotes the (non-normalized) trigonometric system. We let $I$ (respectively, $J$) be the set of nonzero coordinates of the vector $\theta$ (respectively, $\tau$), and denote finally by $\theta_{I}$ (respectively,  $\tau_{J}$) the vector of dimension $|I|$ (respectively, $|J|$) which
contains the nonzero coordinates of $\theta$ (respectively, $|\tau|$). Recall that all densities are defined with respect to the prior $\pi$, which is made explicit in Subsection \ref{PLUS}.
\medskip

For a generic $h\in\mathcal{F}_{m_{h}}(C+1)$, given $\tau\in\mathcal{S}_{1,+}^{p}$ and $s>0$, we let the density $\mbox{dens}_{s}(h|\tau,m_{h})$ with respect to $\pi$ be defined by
\begin{align*}
&\mbox{dens}_{s}(h|\tau,m_{h}) \\
&\quad \propto \exp\left[-\frac{1}{2s^2}
          \sum_{j=1}^{m_{h}}\left(\beta_{h,j}-\tilde{\beta}_{j}(\tau,m_{h}) \right)^2
\right]
\mathbf{1}{\left[\sum_{j=1}^{m_h}j |\beta_{h,j}| \leq C+1 \right]},
\end{align*}
where the $\tilde{\beta}_{j}(\tau,m_{h})$ are the empirical least square coefficients
given by
$$ \left\{\tilde{\beta}_{j}(\tau,m_{h})\right\}_{j=1, \hdots, m_h}
\in  \arg\min_{b\in \mathbb R^{m_h}} \sum_{i=1}^{n} \left(Y_i - \sum_{j=1}^{m_{h}} b_{j}
\varphi_{j}(\tau^{T}\mathbf{X}_{i}) \right)^{2} .$$
In the experiments, we fixed $s=0.1$. Note that simulating with respect to $\mbox{dens}_{s}(h|\tau,m_{h})$ is an easy task, since one just needs to compute a least square estimate and then draw from a truncated Gaussian distribution.
\subsection{Description of $k_{1}$}
We take
\begin{align*}
k_{1}\left(\cdot|(\theta,f)\right) & =
\frac{2k_{1,=}\left(\cdot|(\theta,f)\right)+k_{1,+}\left(\cdot|(\theta,f)\right)}
{3}\,\mathbf{1}_{[|I|=1]}
\\
 & \quad + \frac{k_{1,-}\left(\cdot|(\theta,f)\right)+2k_{1,=}\left(\cdot|(\theta,f)\right)
+k_{1,+}\left(\cdot|(\theta,f)\right)}{4}\,\mathbf{1}_{[1<|I|<p]}
\\
  &   \quad + \frac{k_{1,-}\left(\cdot|(\theta,f)\right)+2k_{1,=}\left(\cdot|(\theta,f)\right)}
{3}\,\mathbf{1}_{[|I|=p]}.
\end{align*}
Roughly, the idea is that $k_{1,-}$ tries to remove one component in $\theta$,
$k_{1,=}$ keeps the same number of components, whereas $k_{1,+}$ adds one component. The density $k_{1,=}$ takes the form
$$ k_{1,=}\left((\tau,h)|(\theta,f)\right)  =  k_{1,=}(\tau|\theta) {\rm dens}_{s}(h|\tau,m_{f}) .$$
The density $ k_{1,=}(.|\theta)$ is the density of $\tau$ when $J=I$ and
$$ \tau_{I} = \frac{\theta_{I}+E}{\|\theta_{I}+E\|_{1}} {\rm
sgn}\left((\theta_{I}+E)_{j(\theta_I+E)}\right), $$
where $E =(E_{1},\hdots,E_{|I|})$ and the $E_{i}$ are independent random variables uniformly distributed in $[-\delta,\delta]$. Throughout, the value of $\delta$ was fixed at $0.5$. It is noteworthy that when we change the parameter from $\theta$ to $\tau$, then we also change the function from $f$ to $h$. Thus, with this procedure, the link function $h$ is more ``adapted'' to $\tau$ and the subsequent move is more likely to be accepted in the Hastings-Metropolis algorithm.
\medskip

In the case where we are to remove one component, $k_{1,-}$ is given by
$$ k_{1,-}\left((\tau,h) |(\theta,f)\right) = \sum_{j\in I} c_{j}
\mathbf{1}_{[\tau=\theta_{-j}]}
                                          {\rm dens}_{s}(h|\tau,m_{f}),$$
where $\theta_{-j}$ is just obtained from $\theta$ by setting the
$j$-th component to 0 and by renormalizing the parameter in order to have $\|\theta_{-j}\|_{1}=1$. We set
$$ c_{j} = \frac{\exp\left(-|\theta_{j}|\right) \mathbf{1}_{[|\theta_{j}|<\delta]}
}{\sum_{\ell\in I} \exp\left(-|\theta_{\ell}|\right)
                    \mathbf{1}_{[|\theta_{\ell}|<\delta]}  } .$$
The idea is that smaller components are more likely to be removed than larger ones. Finally, the density $k_{1,+}$  takes the form
$$ k_{1,+}\left((\tau,h) |(\theta,f)\right) = \sum_{j\notin I} c'_{j} 
                                       \mathbf{1}_{[\tau_{-j}=\theta]}                           
\frac{\mathbf{1}_{[|\tau_{j}|<\delta]}}{2\delta} {\rm dens}_{s}(h|\tau,m_{f}).$$
We set
$$ c'_{j} = \frac{\exp\left(\left|\sum_{i=1}^{n}\left(Y_{i}-f(\theta^{T}\bX_{i})\right)
(\bX_{i})_{j}\right|\right) }
   {\sum_{\ell\notin I} \exp\left(\left|\sum_{i=1}^{n}\left(Y_{i}-f(\theta^{T}\bX_{i})\right)
(\bX_{i})_{\ell}\right|\right) } $$
where $(\bX_{i})_{j}$ denotes the $j$-th component of $\mathbf{X}_{i}$. In words, the idea is that a new nonzero coordinate in $\theta$ is more likely to be interesting in the model if the corresponding feature is correlated with the current residual.
\subsection{Description of $k_{2}$}
In the same spirit, we let the conditional density $k_2$ be defined by
\begin{align*}
k_{2}\left(\cdot|(\theta,f)\right) & =
\frac{2k_{2,=}\left(\cdot|(\theta,f)\right)+k_{2,+}\left(\cdot|(\theta,f)\right)}
{3}\,\mathbf{1}_{[m_{f}=1]}
\\
 & \quad + \frac{k_{2,-}\left(\cdot|(\theta,f)\right)+2k_{2,=}\left(\cdot|(\theta,f)\right)
+k_{2,+}\left(\cdot|(\theta,f)\right)}{4}\,\mathbf{1}_{[1<m_{f}<n]}
\\
  &   \quad + \frac{k_{2,-}\left(\cdot|(\theta,f)\right)+2k_{2,=}\left(\cdot|(\theta,f)\right)}
{3}\,\mathbf{1}_{[m_{f}=n]}.
\end{align*}
We choose 
$$k_{2,=}\left((\tau,h)|(\theta,f)\right) = \mathbf{1}_{[\tau=\theta]}
{\rm dens}_{s}(h|\tau,m_{f})$$ 
and
$$ k_{2,+}\left((\tau,h)|(\theta,f)\right) = \mathbf{1}_{[\tau=\theta]}
 {\rm dens}_{s}(h|\tau,m_{f}+1).$$
With this choice, $m_h=m_f+1$, which means that the proposal density tries to add one coefficient in the expansion of $h$, while leaving $\theta$ unchanged. Finally
$$ k_{2,-}\left((\tau,h)|(\theta,f)\right) = \mathbf{1}_{[\tau=\theta]}
 {\rm dens}_{s}(h|\tau,m_{f}-1),$$ 
 and the proposal tries to remove one coefficient in $h$.
\paragraph{Acknowledgments.} The authors thank three referees for valuable comments and insightful suggestions, which lead to a substantial improvement of the paper. They also thank John O'Quigley for his careful reading of the article. 
\bibliographystyle{plain}
\bibliography{biblio-singleindex}
\end{document}